\documentclass[centertags,leqno]{article}

\usepackage{amsmath}
\usepackage{amssymb}
\usepackage{latexsym}
\usepackage{amsxtra}
\usepackage{amscd}
\usepackage{theorem}

\usepackage{curves}
\usepackage{epic}

\theoremstyle{plain}

{\theorembodyfont{\slshape}

        \newtheorem{thm}{Theorem}[section]
        \newtheorem{cor}[thm]{Corollary}
        \newtheorem{lemma}[thm]{Lemma}
        \newtheorem{prop}[thm]{Proposition}
        
         \newtheorem{con}[thm]{Conjecture} 
        
}

{\theorembodyfont{\rmfamily}

        \newtheorem{defn}[thm]{Definition}
        
        \newtheorem{rem}[thm]{Remark}

}

\renewcommand{\em}{\sl}

\newcommand{\proof}{\noindent {\bf Proof:\ }}
\newcommand{\Endproof}{\hspace*{\fill} $\Box$ \vspace{1ex} \noindent }

\def\ha{{{\gamma} }}
\def\hb{{{\delta} }}

\def\e{{\rm e}}

\def\rk{{\rm rk}}
\def\F{{\cal F}}
\def\O{{\cal O}}
\def\B{{\cal B}}

\newcommand{\Res}{ {\rm Res} }
\newcommand{\LocSys}{ {\rm LocSys} }

\newcommand{\Mod}{ {\rm Mod} }

\newcommand{\p}{ {\frak p} }

\def\K{{\cal K}}
\newcommand{\CC}{{\Bbb C}}

\def\FF{{\cal F}}
 \def\PP{{\Bbb P}}
\def\QQ{{\Bbb Q}}
\def\ZZ{{\Bbb Z}}

\def\GL{{\rm GL}}

\def\im{{\rm im}}

\def\Mon{{\rm Mon}}
\def\S{{\cal S}}

\def\NN{{\Bbb N}}

\def\L{{\cal L}}
\def\G{{\cal G}}
\def\A{{\bf A}}
\def\a{{\bf a}}

\def\K{{\cal K}}

\newcommand{\tr}{{\mbox{\rm {\small tr}}}}

\title{On the middle convolution}

\author{
   Michael Dettweiler \thanks{The first 
 author gratefully acknowledges 
financial support from the  Deutsche Forschungsgemeinschaft DFG.}\\ 
  {} \and
   Stefan Reiter\thanks{The  second author gratefully acknowledges 
financial support from the  Research Training
Network
(Galois Theory and Explicit Methods in Arithmetic) of  the European
Community.}\\{}}

\begin{document}

\maketitle

\begin{abstract} In \cite{DR}, a purely
algebraic analogon 
 of Katz' middle convolution functor 
 (see \cite{Katz97}) is 
given. It is denoted by $MC_\lambda$.
In this paper, we present a cohomological
interpretation of $MC_\lambda$   and find an 
explicit Riemann-Hilbert correspondence for this functor.
   This leads to 
an algorithm for the  construction 
of
 Fuchsian systems corresponding to irreducible rigid 
local systems under the Riemann-Hilbert correspondence.
Also, we  describe
the effect of $MC_\lambda$ on the $p$-curvatures
and find new examples of differential equations for which 
the Grothendieck-Katz $p$-curvature conjecture holds.
\end{abstract}

\section{Introduction}

Let $D$ be a complex ordinary differential equation of 
order $n$ or, equivalently,
a linear system of differential equations of rank $n.$
 Let $T=\{t_1,\ldots,t_r\}\subseteq \CC$
denote  the set of finite singularities of $D$ and let 
$\gamma$ be a closed path in $X=\CC\setminus T.$
Analytic  continuation of  a fundamental 
matrix $F$ of $D$ along $\gamma$ transforms $F$ into $F\cdot A,$
with $A\in \GL_n(\CC)$ uniquely determined. One calls $A$ the {\em 
monodromy 
matrix} of $D$ with respect to $\gamma.$ In other words, if 
$\FF$ denotes the  local system on $X,$
formed by the solutions of $D,$
then
$A$
describes the monodromy of  $\FF$  along $\gamma.$

Since Riemann's investigations on the 
hypergeometric equations (\cite{Riemann57}), the use of monodromy 
is one of the most powerful
tools in the investigation of integrable differential equations.\\

 It is a basic fact, already used by Riemann,
that the solutions of the hypergeometric 
differential equations  give rise 
to a physically 
rigid local system, see \cite{Katz97}, Introduction. 
This means essentially, that the global behaviour
of the solutions under analytic continuation is determined 
by the local behaviour at the singularities (including $\infty$).

A description of
all  irreducible and physically 
rigid local systems on the 
punctured affine line
was given by
Katz \cite{Katz97}. The main tool herefore 
is a middle 
convolution functor on the category
of  perverse sheaves, loc. cit., Chap. 5.
This functor is denoted by 
$MC_\chi,$ for $\chi$ a one-dimensional
representation of $\pi_1({\Bbb G}_m).$ It
  preserves
important properties of local systems 
like the index of rigidity and 
irreducibility,  but in general, $MC_\chi$
changes the rank and the monodromy group.
As an application, Katz shows that
any irreducible rigid local system on the 
punctured affine line can be obtained from a one-dimensional 
local system by applying iteratively  a suitable sequence of middle
convolutions $MC_{\chi_i}$ and scalar multiplications, loc. cit., Chap. 6.
Since the effect of $MC_{\chi_i}$ on the local monodromy can be determined 
via Laumon's theory of $l$-adic Fourier transform, this leads 
to an existence  algorithm for rigid local systems, loc. cit., Section
6.4.\\

In \cite{DR}, the authors give a a purely algebraic analogon
 of the functor $MC_\chi$ 
 (the construction is reviewed in Section \ref{convo}).
This analogous functor 
is a functor of the category $\Mod(K[F_r])$ of
modules of 
the free group $F_r$ on $r$ generators to itself. 
It  depends on 
a scalar $\lambda\in \CC^\times$ and is denoted by $MC_\lambda.$  

One has 
the equivalence between $\Mod(\CC[F_r])\cong 
\Mod(\CC[\pi_1(X,x_0)])$
and the category 
$\LocSys(X)$ of local systems on $X,$ see Section  \ref{localsyst}.
Then,
$MC_\lambda$ translates into a functor
of the category of local systems on the $r$-punctured 
complex affine line $X$ to itself, sending 
a local system $\F$ to $MC_\lambda(\F),$ 
see Section \ref{42}.
It follows from the results of 
\cite{DR} that $MC_\lambda$ (viewed as a functor 
on the category of local systems on $X$) has 
 analogous properties as Katz' functor $MC_\chi,$
where $\chi$ is the representation, sending 
the standard  generator of $\pi_1({\Bbb G}_m(\CC))$ to $\lambda.$ 
This leads 
to a new and elementary proof of Katz' 
existence algorithm for rigid local systems,
see \cite{DR}, Chap. 4. Similar results are
obtained in \cite{Voelklein01}.\\

 It is the aim of this paper to give answers to 
the following   problems:\\

\noindent  {\bf Problem 1:}  Give a cohomological 
interpretation of $MC_\lambda(\F),$ explaining 
the formal similarity between 
$MC_\lambda$ and
Katz' functor $MC_\chi.$\\

By the work of Deligne, it is known that 
the category of complex local systems on $X=\CC\setminus T$ is equivalent 
to the category of ordinary complex differential equations 
 with 
polynomial coefficients having only regular singularities at 
the missing points (including $\infty$) and no singularities in $X,$
 see 
\cite{Deligne70}. This equivalence is called the {\em Riemann-Hilbert 
correspondence}. This leads to\\

 \noindent \noindent  {\bf Problem 2:} 
  Given a differential system having only regular singularities
and whose local system formed by its solutions is $\FF,$
 find a differential system having only regular singularities
such that the local system of its solutions  
is $MC_\lambda(\F).$\\

As it turns out  both problems are closely related to the 
 cohomology of the  locally trivial 
fibration ${\rm p}_2:E\to X,$ defined below. The first one is 
related to the singular cohomology and second one 
to the de Rham cohomology of ${\rm p}_2.$ \\

Motivated by Katz' 
description of $MC_\chi(\FF)$ in \cite{Katz97}, Chap. 5.1,
 we give a solution 
to Problem (1) in Section \ref{s2}. For this, let $X=\CC\setminus T$
 be as above,
$$E=\{(x,y)\in \CC^2 \mid x,y\not= t_i,\, i=1,\ldots,r,\, x\not= y\},$$
 ${\rm p}_i:E\to X,\, i=1,2,$ be the $i$-th projection, $$q:E\to \CC^\times,\,
(x,y)\mapsto y-x,$$ $j: E \to \PP^1(\CC)\times X$ the tautological 
inclusion
 and 
$\bar{p}_2: \PP^1(\CC)\times X \to X$ the  (second) projection onto $X.$ 
Moreover, let  $\L_\lambda$ denote the  Kummer sheaf
associated to the representation, which sends 
a generator of $\pi_1(\CC^\times)$ to $\lambda$ 
(see Definition \ref{Kummer}).  
The following theorem is proved in Section \ref{cc}, using singular 
sheaf cohomology (see 
Theorem \ref{ccc}):

\begin{thm}\label{ka} Let $\F$ be a 
local system on $X$ and $\lambda \in \CC^\times\setminus 1.$
Then 
$ MC_\lambda(\F)$ is isomorphic to the higher 
direct image sheaf 
$ R^1 (\bar{p}_2)_*(j_*({\rm p}_1^*(\F)\otimes 
q^*(\L_\lambda))).$
 \end{thm}
 
The idea of the proof is to relate the construction of 
$MC_\lambda$
to the 
 group (resp. singular) cohomology of the locally trivial 
fibration ${\rm p}_2:E\to X,$  where one can explicitly 
work with crossed homomorphisms. We do not use the standard base
but a twisted base which arises from the use of 
commutators, also called Pochhammer contours.
The Pochhammer contours are crucial 
in the further investigation of the convolution 
in terms of Fuchsian systems (see the Theorem \ref{1} below and 
Remark \ref{contrast}).
Translating Theorem \ref{ka} into the language of
perverse sheaves, one rediscovers Katz' original 
construction, see \cite{Katz97}, 5.1.7.\\

In Section \ref{Eule}
we consider Problem (2): In 
\cite{DR}, Appendix A, an additive version of Katz' functor is defined. It 
depends on a scalar $\mu \in \CC$ and is 
denoted by $mc_\mu.$ By definition,  $mc_\mu$ is  nothing 
else then a transformation of tuples of matrices 
$$(a_1,\ldots,a_r) \in
(\CC^{n\times n})^r \mapsto mc_\mu(a_1,\ldots,a_r)\in (\CC^{m\times m})^r.$$

Any choice of elements $t_1,\ldots,t_r\in \CC,$ 
together with a tuple of matrices $\a:=(a_1,\ldots,a_r)\in (\CC^{n\times n})^r,$
yields a Fuchsian system 
$$D_\a: Y'=\sum_{i=1}^r \frac{a_i}{x-t_i}Y.$$ 
Then, $mc_\mu$   translates
 into a transformation of   Fuchsian systems,
sending $D_\a$ to $D_{mc_\mu(\a)}.$  This transformation 
will be called the 
{\em middle convolution} of Fuchsian systems.
The
tuple of monodromy generators of $D_\a$ will be  denoted by 
$\Mon(D_\a)$ (see Section \ref{fuchs}). One obtains the following result,
see Theorem \ref{monodromy}:

\begin{thm}\label{1} {\rm (Riemann-Hilbert correspondence for $MC_\lambda$)}
Let $\mu \in \CC\setminus \ZZ,$ $\lambda =
\e^{2\pi i \mu}$ and 
$\a:=(a_1,\ldots, a_r),\, a_i\in \CC^{n\times n},$
such that  $\Mon(D_\a)
=(A_1,\ldots,A_r)\in \GL_n(\CC)^r.$ Assume that
$$\rk(a_i) = \rk(A_i-1),\,\,\,\,
\rk(a_1+\cdots + a_r+\mu)=\rk( \lambda\cdot A_1\cdots A_r-1)$$ and that
$\langle A_1,\ldots,A_r \rangle$ generates an irreducible subgroup
of $\GL_n(\CC)$ such that at least two elements $A_i$ are $\not=1.$ Then
$$    \Mon(D_{mc_{\mu-1}(\a)})=MC_{\lambda}(\Mon(D_\a)).$$
\end{thm}

Let $\FF_\a$ (resp. $\F_{mc_{\mu-1}(\a)}$)
denote the local system, formed by the solutions 
of $D_\a$ (resp. $D_{mc_{\mu-1}(\a)}$).
Since $\Mon(D_\a)$ (resp. $\Mon(D_{mc_{\mu-1}(\a)})$)
describes the monodromy of the 
local system $\FF_\a$ (resp. $\F_{mc_{\mu-1}(\a)}$) (see Remark 
\ref{rightleft}),
Theorem \ref{1} yields 
$$ \F_{mc_{\mu-1}(\a)}\cong MC_\lambda(\F_\a).$$
Thus we have obtained the 
 Riemann-Hilbert correspondence (Problem 2)
 under the assumptions of Theorem \ref{1}.
These assumptions are rather mild and can be further
weakened (see the remark following Theorem \ref{monodromy}). 

The main idea of the proof of Theorem \ref{1} is to use
Euler transformations, in order to  construct a  suitable
period matrix $I^\mu,$ describing the pairing 
between the homology and the de Rham cohomology with 
coefficients in $({\rm p}_1^*(\F_\a)\otimes 
q^*(\L_\lambda))|_{X(y_0)},$ resp. 
$({\rm p}_1^*(\F_\a)\otimes 
q^*(\L_\lambda))^\vee|_{X(y_0)},$ where 
$X(y_0)$ denotes the fibre of ${\rm p}_2:E\to X$ 
over $y_0.$   The  columns of 
$I^\mu$ are solutions of a Fuchsian system 
$D_{c_\mu(\a)}$ (called the {\em convolution} of $D_\a$ with $\mu$) such that 
the middle convolution
$D_{mc_\mu(\a)}$ is a factor system of $D_{c_\mu(\a)}.$
The rows of $I^\mu$
have an interpretation in 
terms of crossed homomorphisms (already used in 
Section \ref{s2})   
which makes it possible to compute the monodromy
of $D_{c_\mu(\a)}$ and $D_{mc_\mu(\a)}.$ \\

Finally, we give some applications of our 
methods (Section \ref{apl}): 
From Theorem \ref{1},
one obtains an algorithm for the  construction of  Fuchsian 
systems corresponding to irreducible rigid local systems
under the Riemann-Hilbert correspondence, see
Section \ref{Rigidesyst}.
As a byproduct,
one obtains integral expressions 
for the solutions of these Fuchsian systems.  Compare to
the work of 
Haraoka and Yokoyama (\cite{HY}, \cite{Yokoyama02}) who use 
 a different approach 
(in the case of semisimple monodromy) to obtain 
integral expression of such solutions.

Then
 we apply Theorem \ref{1} to 
the construction problem
of  differential systems which arise from
geometry: These are  
differential systems which  arise from iterated extensions 
of subfactors of Gau\ss-Manin connections 
 (see Section \ref{Kurven} for the 
definition). These  differential 
 systems have many favorable properties. For example,
under some additional assumptions (the connectivity 
of motivic Galois groups), such a system satisfies the 
Grothendieck-Katz $p$-curvature 
conjecture which makes it possible
to construct the Lie algebra of its differential Galois group
from its $p$-curvatures, see Andr\'e \cite{Andre02}, Theorem 0.7.1.

Using  results of Andr\'e \cite{Andre89}, one obtains the 
following result (Theorem \ref{coming}):

\begin{thm} Let $K$ be a number field,
 $\a=(a_1,\ldots,a_r),\, a_i\in
 K^{n\times n},\,\mu \in \QQ,$ such that the conditions
of Theorem \ref{1} hold for $D_\a.$ 
If $D_\a$ is arising from geometry, then
$D_{mc_\mu(\a)}$ is arising from geometry. 
\end{thm}

This makes it  possible to construct explicitly 
a large  number of differential systems which arise from 
geometry. One could start from any 
 differential system with finite monodromy 
(which automatically arises from geometry) and apply the 
convolution $mc_\mu,\, \mu \in \QQ,$ to it. 
In Section \ref{Kurven}, we consider 
examples which are derived from   Lam\'e equations with 
finite monodromy, related to the work of Baldassari 
\cite{Baldassari80} and 
Beukers and van der Waall \cite{BW}, \cite{vdW}. 
This leads to new (non-rigid) examples of differential systems
for which the Grothendieck-Katz $p$-curvature conjecture is true,
see Corollary \ref{ak}. See Katz \cite{Katz97}, Chap. 9, for a proof
of the Grothendieck $p$-curvature 
conjecture for Fuchsian systems corresponding
to irreducible rigid local systems. 

As another application,
 we investigate the effect of the convolution on the 
$p$-curvature ($p$ a prime) of a Fuchsian system defined over $\QQ.$ 
One obtains a simple formula for the computation of the 
$p$-curvature matrices  (Lemma \ref{recu}) 
and the following result, see Section \ref{pcurv} for 
definitions and Theorem \ref{nilpo}:

\begin{thm}\label{nilp} Let $K$ be a number field 
and $\p$ a prime of $K$ lying over p.
Let $\mu \in \QQ$ and  
 $\a=(a_1,\ldots,a_r),\, a_i \in K^{n\times n},$ such that 
the $p$-curvature matrix $a{(\p)}$ of $D_\a$ satisfies 
$a{(\p)}^k=0.$ 
Then the $p$-curvature matrix $mc_\mu(a{(\p)})$
of the 
convoluted Fuchsian system  $D_{mc_\mu(\a)}$
satisfies $mc_\mu(a{(\p)})^{k+2}=0.$
\end{thm}

The
 crucial observation here is, that the convolution 
   $D_{c_\mu(\a)}$ is a differential system in 
Okubo normal form (see Section \ref{pcurv} for definition). For these systems
there exists a closed formula for the computation of
the $p$-curvature matrices (Lemma \ref{recu}).
Theorem \ref{nilp} is interesting in view of the Bombieri-Dwork
conjecture which relates the nilpotence of the 
$p$-curvatures to the geometric nature of a differential 
equation, see Section \ref{pcurv}. Also,
information on the $p$-adic radius of solvability is encoded 
in the nilpotence degree of the $p$-curvatures, see \cite{Andre89}, Chap.~4.\\

The authors thank 
Y. Andr\'e and D. Bertrand for valuable conversations
and 
P. Deligne for suggesting the use of crossed homomorphisms
 for a geometric 
interpretation of $MC_\lambda.$

The second author wants to thank the Department of 
Mathematics of the University of Paris 6 (Jussieu)
  - especially  Y. Andr´e and D. Bertrand - and  
the department of Mathematics of the University
 Toulouse II (G.R.I.M.M.) - especially the group ALGO -  for their
hospitality.

\section{Definition and properties
of the middle convolution functor $MC_\lambda$}\label{convo}
In this section, we recall the algebraic construction 
of  the  multiplicative version 
of the convolution functor defined 
in \cite{DR}. We actually consider a slight 
modification of the multiplicative version of the convolution.
This modification  
is just of formal nature and due to the topological 
setup used in the later 
sections.\\

We will use the following notations and conventions throughout the 
paper:
Let $K$ be 
 a field and $G$ a group. The category 
of finite dimensional 
left-$G$-modules is denoted by $\Mod(K[G]).$ Mostly, we do not 
distinguish notationally between an  element of $\Mod(K[G])$
and its underlying vector space.
Let $V$ be an element $\Mod(K[G])$ corresponding to 
a representation $\rho: G\to \GL(V)$
 and $W$ a $K$ vector space such 
that one has a perfect pairing 
$$<,>:W\times V \to K. $$ 
Then $W$  turns into 
a $G$-module, where $g$ acts via 
the unique linear 
transformation  $\rho^\vee(g)$ such that 
$<\rho^\vee(g)w,\rho(g)v>=<w,v>$ for all $ w\in W $ and
all $v\in V.$ 
We refer to $W$ as the dual module of $V$ with respect to 
$<,>$ and often  denote it  $V^\vee.$
If for $g\in G,$ the linear transformation
$\rho(g)$ is  a given element 
$A\in \GL(V),$ then 
 we write $A^\vee\in \GL(W)$ for $\rho^\vee(g).$\\

\subsection{Definition of $MC_\lambda$}\label{katz}

Let  $F_r$ denote  the free group on $r$ generators $f_1,\ldots,f_r.$ 
 An element 
in $\Mod(K[F_r])$ is viewed as a pair $(\A,V),$ 
where $V$ is a vector space over $K$ and $\A=(A_1,\ldots,A_r)$ is 
an element of 
$\GL(V)^r$ such that $f_i$ acts on $V$ via $A_i,\, i=1,\ldots,r.$
For $(\A,V)\in \Mod(K[F_r]),$ where $\A=(A_1,\ldots,A_r)\in 
\GL(V)^r,$ and $\lambda \in K^\times$ one can 
construct an element $(C_\lambda(\A),V^r)\in \Mod(K[F_r]),$
$C_\lambda(\A)=(B_1,\ldots,B_r)\in \GL(V^r)^r,$ as follows:
For $k=1,\ldots,r,$ 
$B_k$ maps a vector $(v_1,\ldots,v_r)^{\rm tr}$ $\in V^r$
to
\[ \left( \begin{array}{ccccccccc}
                  1 & 0 &  & \ldots& & 0\\
                   & \ddots &  & & &\\
                    & & 1 &&&\\
               \lambda (A_1-1) & \ldots& \lambda (A_{k-1}-1)  & \lambda A_{k} & (A_{k+1}-1) & \ldots 
&   (A_r-1) \\
     &&&&1&&\\
 %            0 & \ldots & & 1 & 0 &   \ldots  \\     
               &   &  & && \ddots  &   \\             
                   0 &  &  & \ldots& &0 & 1
          \end{array} \right)\left(\begin{array}{c}
v_1\\
\vdots\\
\\\vdots\\
\\\vdots\\
\\ v_r\end{array}\right)
.\]

 We set $C_\lambda(\A):=(B_1,\ldots,B_r).$
There are the following $\langle B_1,\ldots,B_r \rangle$-invariant 
subspaces 
of  $V^r:$

\[ \K_k = \left( \begin{array}{c}
          0 \\
          \vdots \\
          0 \\
          \ker(A_k-1) \\
            0\\
           \vdots \\
           0 
        \end{array} \right)  \quad \mbox{({\it k}-th entry)},\, k=1,\dots,r,\]
and 
\[     \L=\cap_{k=1}^r \ker (B_k-1)={\rm ker}(B_1\cdots B_r - 1).
\]
Let $\K:=\oplus_{i=1}^r\K_i.$

If $\lambda \not=1,$ then 
 $$\L=
\langle \left( \begin{array}{c}
            A_2 \cdots A_{r} v \\
                A_3 \cdots A_{r}v \\
               \vdots \\
                        v
            \end{array} \right) \mid v \in \ker(\lambda\cdot A_1\cdots A_r-1) \rangle.$$ 
and $$\K+ \L=\K \oplus \L.$$

\begin{defn}{\rm Let $V=(\A,V)\in \Mod(K[F_r]).$

i) We call the $K[F_r]$-module
$ C_{\lambda}(V):=(C_\lambda(\A),V^r)$
the {\em convolution} of $V$ with $\lambda.$

 ii) Let 
$MC_\lambda(\A):=(\tilde{B}_1,\dots,\tilde{B}_r)\in 
\GL(V^r/(\K+\L))^r,$ where $\tilde{B}_k$ is induced by the 
action of $B_k$ on $V^r/(\K+\L).$
The $K[F_r]$-module
$MC_{\lambda}(V):=(MC_\lambda(\A),V^r/(\K+\L))$ 
 is called the
{\em  middle convolution} of $(A_1,\ldots,A_r)$ with 
$\lambda.$}\end{defn}

\noindent{\bf Remark:} In \cite{DR}, we use the same construction, with
the difference that 
the $k$-th block row of $B_k$ is 
$$(  (A_1-1) , \ldots,  (A_{k-1}-1)  , \lambda A_{k} ,\lambda
(A_{k+1}-1) , \ldots 
, \lambda  (A_r-1)).$$

\subsection{Properties of $MC_\lambda$} \label{property}

Let $V\to V'$ be a morphism of $F_r$-modules. This clearly
induces a morphism $C_\lambda(V)\to C_\lambda(V').$ Since
the subspaces $\K$ and $\L$ of $C_\lambda(V)$ 
are mapped to their corresponding subspaces 
$\K'$ and $\L'$ of $C_\lambda(V')$ this induces a morphism
 $MC_\lambda(V)\to MC_\lambda(V').$ The following proposition
is easy to prove, compare to \cite{DR}, Proposition 2.6 and Lemma 2.8:

\begin{prop}\label{deli2}  Let $\lambda\in K^\times.$
The transformation $V\mapsto MC_\lambda(V)$ 
(resp. $V\mapsto C_\lambda(V)$)
 is a covariant, end-exact,  functor of $\Mod(K[F_r])$
to itself.
\end{prop}

\begin{defn}{\rm Let  $V=(\A,V)\in \Mod(K[F_r]),$
where $\A=(A_1,\ldots,A_r)\in \GL(V)^r.$ 
We say that $V$ satisfies $(*)$ if 
$$ \bigcap_{j\not= i }{\rm ker}(A_j-1) \cap {\rm ker} (\tau A_i -1) =0,\; i=1,\ldots , r
, \;\forall \tau \in K^\times .$$

Let 
 ${\cal U}_i(\tau):=
 \sum_{j \neq i} {\rm im}(A_j-1) +{\rm im} (\tau A_i-1),
\; i=1,\ldots,r, \;\tau \in K^\times. $
We say that  $V$ satisfies $(**)$ if
$$ {\rm dim}({\cal U}_i(\tau))=\dim (V) ,\; 
i=1,\ldots,r, \;\forall \tau \in K^\times. $$
}\end{defn}

\noindent {\bf Remark:} The conditions $(*)$ and $(**)$ say, 
that  $V$
has no $1$-dimensional factors and/or submodules 
with the property 
 that only one (or none)
of the $A_i$ act non-trivially.\\

\begin{thm} \label{eigen} Let $V=(\A,V)\in \Mod(K[F_r]),$ where 
$\A=(A_1,\dots,A_r) \in \GL(V)^r$
and $\lambda \in K^\times.$\\

i) If 
$\lambda \neq 1,$ then $${\rm dim}(MC_\lambda(V))= \sum_{k=1}^{r} \rk (A_k-1)-
 ({\rm dim}(V)-\rk(\lambda\cdot A_1\ldots A_r-1)).$$

ii) If $\lambda_1,\, \lambda_2\in K^\times$ such that $\lambda_1\lambda_2=\lambda$ and 
 $(*)$ and $(**)$ hold for $V,$ then 
    \[ MC_{\lambda_2}MC_{\lambda_1}(V)\cong MC_{\lambda}(V).\]

iii) Under the assumptions of ii), if
$V$ is irreducible, then
$MC_{\lambda}(V)$ is irreducible.\\

iv)   Let   ${\cal B}_r=\langle Q_1,\ldots,Q_{r-1} \rangle $ be 
the abstract Artin braid
group, where  the generators $Q_1, \ldots,Q_{r-1}$ of ${\cal B}_r$ act in the 
following way on tuples $(g_1,\ldots,g_r)\in G^r$
(where $G$ is a group):
\begin{equation}\label{Q_i}
 Q_i(g_1,\ldots,g_r)=(g_1,\ldots,g_{i-1},g_ig_{i+1}g_{i}^{-1},g_i,g_{i+2},
\ldots,g_r),\,\, i=1,\ldots,r-1.\end{equation}
For any $Q \in{\cal B}_r$ there exists a $B\in \GL(V^r/(\K+\L))$ such that
$$MC_\lambda(Q(\A))= Q(MC_\lambda(\A))^B, $$
where $B$ acts via component-wise conjugation.\\

v) Let  $K=\CC,$ $\lambda \in \CC$ 
 be a root of unity  and  
$MC_{\lambda}(\A)=(\tilde{B}_1,\dots,\tilde{B}_r).$ If 
$\langle A_1,\ldots,A_r\rangle$ respects an hermitean form, 
then $\langle \tilde{B}_1,\dots,\tilde{B}_r \rangle $ 
respects an hermitean form.\\

vi) Let the characteristic of $K$ be different
from $2$ and $MC_{-1}(\A)=(\tilde{B}_1,\dots,\tilde{B}_r).$  
If 
$\langle A_1,\ldots,A_r\rangle$ respects an orthogonal
  (resp. symplectic) form,
then $\langle \tilde{B}_1,\dots,\tilde{B}_r \rangle $ 
respects a
   symplectic (resp. orthogonal) form.
\end{thm}

\proof i)-iv) follow
 analogously to \cite{DR}, Lemma 2.7, Lemma A.4, Theorem 3.5, Corollary 3.6 and 
Theorem 5.1 (in this order). The claims v) and 
vi) follow from Lemma \ref{form2} below. 
\Endproof

\noindent {\bf Remark:}
 The Jordan canonical  forms of $\tilde{B}_k$ can be computed as in
\cite{Katz97},  Chap. 6 (using \cite{DR}, Lemma 4.1).

\begin{lemma}\label{form2}  Let $\A=(A_1,\ldots,A_r),\, A_k \in \GL_n(K),$
$\lambda \in K^\times$ and  $C_\lambda(\A)=(B_1,\ldots,B_k).$
  Let ${\frak G}$ be an invariant form under $A_i,$ i.e.
  $ A_i^{\tr} {\frak G} A_i={\frak G},\, i=1,\ldots,r.$ 
  Then \[ B_k^{\tr} {\frak H}B_k ={\frak H},\, k=1,\ldots,r, \] 
  where
  \[ {\frak H}_{i,i} = {\frak G} \lambda^{1/2} (A_i^{-1}-1) (A_i-\lambda^{-1}) \]
  and
  \[ {\frak H}_{i,j}= {\frak G}  \lambda^{-1/2} (A_i^{-1}-1) (A_j-1), \quad
{\rm if}\,\,\, i<j ,\]
\[ {\frak H}_{i,j}= {\frak G}  \lambda^{1/2} (A_i^{-1}-1) (A_j-1), \quad{\rm if}\,\,\, i>j. \]
\end{lemma}

\section{The underlying fibration and its cohomology}

We fix a finite 
set  $T:=\{t_1,\ldots,t_r\}\subseteq \CC$ such that  $t_i\not= t_j$ 
for
$i\not= j,$ and set $X:=\CC\setminus T.$
Let $W$ be a topological space  and $I:=[0,1].$
A path in $W$ is a continuous map $\gamma:I\to W.$ If 
$\gamma_1,\gamma_2$ are paths in $W$ such that the endpoint 
of $\gamma_2$ coincides with the initial point of $\gamma_1,$ then 
their product is denoted by $\gamma_1\gamma_2.$ 
If $\gamma$ is a closed path 
in $W$ with initial point $w_0,$ then
$\gamma \in \pi_1(W,w_0)$  will also denote the corresponding 
homotopy class.

\subsection{The underlying fibration}\label{top}

In this subsection we study a fibration whose cohomology will 
lead to the geometric interpretation of $C_\lambda$ and $MC_\lambda$
in Subsection \ref{cc}. The contents of this section are well known,
compare to \cite{Birman74}, Chap. 1, and  
\cite{Voelklein01}. \\

For $n\in \NN,$ consider the configuration space 
$$ \O_{n}:=\{P\subseteq \CC\mid |P|=n\}$$ 
of subsets of $\CC$ of cardinality equal to $r.$ Let further
$$ \O^n:=\{(p_1,\ldots,p_n)\in \CC^n\mid i\not= j \Rightarrow 
p_i\not= p_j\}.$$ 
Since the map $$\O^n\to \O_n,\, (p_1,\ldots,p_n)\mapsto 
\{p_1,\ldots,p_n\},$$ is an unramified covering map (where 
$\O_n$ is equipped with the obvious topology), we will 
consider $\B^n:=\pi_1(\O^n,(b_1,\ldots,b_n))$ as a subgroup
of  $\B_n:=$ $\pi_1(\O^n,$ $\{b_1,\ldots,b_n\})$ via
covering theory.

It is well
known that the fundamental group $\B_n$ is isomorphic
to the abstract Artin braid group, i.e., it 
has a presentation
on  $n-1$ generators  $Q_1,\ldots,Q_{n-1}$ subject to 
the {\it braid relations}
\begin{eqnarray} Q_iQ_j&=&Q_jQ_i\quad \quad \quad \quad {\rm if} \quad |i-j|>1,\nonumber \\
                 Q_iQ_{i+1}Q_i&=&Q_{i+1}Q_iQ_{i+1} \quad {\rm for} \quad
                 i=1,\ldots,n-2.\nonumber \end{eqnarray}
\noindent The group $\B^n$  is isomorphic
to the (abstract) pure Artin braid group and generated by   
 the elements 
$$ Q_{i,j}:=(Q_i^2)^{Q_{i+1}^{-1}\cdots Q_{j-1}^{-1}}=
(Q_{j-1}^2)^{Q_{j-2}\cdots Q_i},$$
where $ 1\leq i<j\leq n.$ \\

Let now $n=r+2.$ We set
$\ha_i:=Q_{1,i+1},\,\, i=1,\ldots,r+1,
$ and 
$\hb_k:=Q_{k+1,r+2},\, k=1,\ldots,r.$
Using the braid relations, one can easily see that for 
$k=1,\ldots,r,$ the following formula holds:
\begin{eqnarray}\nonumber
 ({}^{\hb_k} \ha_1,\ldots,{}^{\hb_k}\ha_{r+1})=
(\ha_1,\ldots, \ha_{k-1},\ha_k^{\ha_{r+1}},
\ha_{k+1}^{[\ha_k,\ha_{r+1}]},\ldots,\ha_r^{[\ha_k,\ha_{r+1}]}, \ha_{r+1}^{\ha_k
\ha_{r+1}})\,,\end{eqnarray}
where ${[\ha_k,\ha_{r+1}]}=
 \ha_k^{-1}\ha_{r+1}^{-1}\ha_k\ha_{r+1}$ and 
${}^{\hb_k}\ha_1=\hb_k\ha\hb_k^{-1},$ see also 
\cite{Birman74}, 1.8.3, and \cite{Voelklein01}.\\ 

Let $T=\{t_1,\ldots, t_r\}\subseteq \CC,$  $X:=\CC\setminus T$ and 
$$E:= \{ (x,y) \in \CC^2\mid x,y \not= t_i,i=1,\ldots,r,\, x\not= y\}.$$
The second projection ${\rm p}_2: E\to X$ 
is a locally trivial fibration. The fibre over $y$ is denoted by
$X(y_0)$ and is via the first projection
identified with 
$X\setminus \{y_0\}.$ One has a commutative diagram 
\[
\begin{array}{ccccc}
E & \stackrel{{\rm p}_2}{\longrightarrow} &{X}\\
\downarrow &&\downarrow\\
\O^{r+2}&\stackrel{{\rm p}}{\longrightarrow} & \O^{r+1}
\end{array},\]
where ${\rm p}(p_1,\ldots,p_{r+2}):=(p_2,\ldots,p_{r+2})$
and the first (resp. second) 
vertical arrow is given by $(x,y)\mapsto (x,t_1,\ldots, t_r,y)$  
(resp.  $y\mapsto (t_1,\ldots,t_r,y)$).

 The long exact 
sequences of homotopy groups, associated to 
locally trivial fibrations, lead then to a commutative diagram
\begin{equation}\label{lum}\begin{array}{ccccccccccc}
 1& \to& \pi_1(X(y_0),x_0)&\to& \pi_1(E,(x_0,y_0))
&\to& \pi_1(X,y_0)&\to& 1\\
&&\downarrow &&\downarrow&&\downarrow&&\\
1&\to &\pi_1({\cal F}_r,(x_0,t_1,\ldots,y_0)))&\to &\B^{r+2}&\to& \B^{r+1}&\to& 1
\end{array},\nonumber
\end{equation}
where $\B^{r+2}=\pi_1(\O^{r+2},(x_0,t_1,\ldots,t_r,y_0)),$ 
$\B^{r+1} =\pi_1(\O^{r+1},(t_1,\ldots,t_r,y_0))$ and
${\cal F}_r$ denotes the the fibre over $(t_1,\ldots,t_r,y_0).$
It is well known,  that the 
 rows are split exact sequences and the vertical arrows are
injective, see \cite{Birman74}. Moreover, one can check that 
$\pi_1({\cal F}_r,(x_0,t_1,\ldots,t_r,y_0)))$ is generated by 
$\ha_1,\ldots,\ha_{r+1}$ and that the image of $\pi_1(X,y_0)$
in $\B^{r+1}$ is generated by $\hb_1,,\ldots,\hb_{r}.$ \\

We 
define $\alpha_1,\ldots,\alpha_{r+1}\in \pi_1(X(y_0),x_0)$
(resp. $\beta_1,\ldots,\beta_r\in  \pi_1(X,y_0)$)
to be the inverse images of $\ha_1,\ldots,\ha_{r+1}$
(resp. $\hb_1,\ldots,\hb_r$) under the first (resp. third)
 vertical arrow. Thus one deduces that for 
$k=1,\ldots,r$ the following formula holds:
\begin{eqnarray}\label{zopf}
 ({}^{\beta_k}\alpha_1,\ldots,{}^{\beta_k}\alpha_{r+1})&=&
 (\alpha_1,\ldots, \alpha_{k-1},\alpha_k^{\alpha_{r+1}},
\alpha_{k+1}^{[\alpha_k,\alpha_{r+1}]},
 \ldots,\alpha_r^{[\alpha_k,\alpha_{r+1}]}, \alpha_{r+1}^{\alpha_k
\alpha_{r+1}}).\nonumber \end{eqnarray}

\subsection{Group cohomology of the fibration}\label{cohomo}

If $G$ is a group and $\rho\to \GL(V)$ is a representation,
then we define the {\em cohomology of $G$ with values 
in the module $V$} to be $H^1(G,V):=C^1(G,V)/B^1(G,V),$
where 
$$C^1(G,V):=\{(\delta: G\to V)\mid \delta(gg')=
\delta(g')+\rho(g')^{-1}\delta(g),\, \forall g,g' \in G \}$$
is the vector space of {\em crossed homomorphisms} and 
$$ B^1(G,{V})=\{(\delta: G\to {V})\mid \exists
v\in {V},\, \delta(g)=v- \rho(g)^{-1}v,\,
\forall g \in G\} $$ is the subspace of {\em exact}
crossed homomorphisms.\\

Let $V$ be a  $\pi_1(X,x_0)$-module, where
$\alpha_i$ acts via $A_i\in \GL(V),\, i=1,\ldots,r.$
Let 
 $\lambda \in \CC,$ 
$\Pi:=\pi_1(X(y_0),x_0)=
\langle \alpha_1,\ldots,\alpha_{r+1}\rangle,$ and 
${V}_\lambda$ be the $\Pi$-module, whose 
underlying vectorspace is $V$ and where $\alpha_1,\ldots,\alpha_r$ 
act via $A_i \in\GL(V)$ and $\alpha_{r+1}$ acts via $\lambda.$
The underlying representation is denoted by $\rho_\lambda.$

\begin{defn}{\rm  The linear map
$$ \tau :C^1(\Pi,{V_\lambda})\to ({V_\lambda})^r,\, \delta \mapsto 
(\delta([\alpha_{r+1},\alpha_1]),\ldots,
\delta([\alpha_{r+1},\alpha_r]))^{\tr}$$
is called the {\em twisted evaluation map.}}\end{defn}

\begin{lemma}\label{3.4}  If $\lambda\not= 1,$
then the kernel of the twisted evaluation map 
$\tau: C^1(\Pi,V_\lambda)\to V_\lambda^{r}$
is $B^1(\Pi,{V_\lambda}).$
\end{lemma}

 \proof
The crossed homomorphism relation implies
$$\delta([\alpha_{r+1},\alpha_i])= (1-\lambda^{-1})   \delta(\alpha_i)
-(1-A_i^{-1})\delta(\alpha_{r+1}).$$ 
So, if  
$\delta([\alpha_{r+1},\alpha_i])=0$ for $i=1,\ldots, r,$ then 
$$ \delta(\alpha_i) =\frac{1}{1-\lambda^{-1}} (1-A_i^{-1})  \delta(\alpha_{r+1})
,\, i=1,\ldots,r.$$ An easy induction 
shows that
$$\delta(\gamma)=\frac{1}{1-\lambda^{-1}}(1-\rho_\lambda(\gamma)^{-1})  \delta(\alpha_{r+1}),$$
so $\delta$ is exact. On the other hand, any vector in ${V_\lambda}$ 
occurs as $\delta(\alpha_{r+1})$ for some $\delta\in C^1(\Pi,{V_\lambda}).$
Therefore, the claim follows from dimension reasons. \Endproof\\

Since 
\begin{eqnarray}
 ({}^{\beta_k}\alpha_1,\ldots,{}^{\beta_k}\alpha_{r+1})&=&
 (\alpha_1,\ldots, \alpha_{k-1},\alpha_k^{\alpha_{r+1}},
\alpha_{k+1}^{[\alpha_k,\alpha_{r+1}]},
 \ldots,\alpha_r^{[\alpha_k,\alpha_{r+1}]}, \alpha_{r+1}^{\alpha_k
\alpha_{r+1}})\nonumber \end{eqnarray}
and  $\rho_\lambda(\alpha_{r+1})=\lambda,$ 
the map which sends
$  \delta $ to 
$\delta \circ \beta_k^{-1}$
is contained in 
$\GL(C^1(\Pi,{V_\lambda})).$ 
Thus, by Lemma \ref{3.4},the association
$$ \beta[\delta]:=[\delta\circ \beta^{-1}]$$
imposes the structure of a $\pi_1(X,y_0)$-module
on $H^1(\Pi,{V_\lambda})$ and, by the same arguments, on 
$H^1(\Pi,{V_\lambda^\vee}).$\\

Consider the pairing 
$$ (V_\lambda^\vee)^r\times V_\lambda^r\to \CC , \,\, ((w_1,\ldots,w_r),(v_1,\ldots,v_r)^\tr)
\mapsto <w_1,v_1>+\cdots + <w_r,v_r>.$$
Let $C_\lambda(\A)=(B_1,\ldots,B_r)$ and $\tilde{C}_\lambda(V_\lambda)$ the 
$\pi_1(X,y_0)$-module whose underlying vector space is $V_\lambda^r,$ on
which $\beta_k$ acts via $B_k.$ 
Let further $\tilde{C}_\lambda(V_\lambda)^\vee$
denote the dual module with respect to the above pairing.

\begin{thm}\label{Bk} The linear map
$$ H^1(\Pi,V_\lambda^\vee) \to \tilde{C}_\lambda(V_\lambda)^\vee,\,\,
[\delta] \mapsto \tau(\delta)$$
is an isomorphism of  $\pi_1(X,y_0)$-modules.
 \end{thm}

\proof  It suffices to show that 
\begin{eqnarray}\label{3.2} \tau(\beta_k^{-1}(\delta))&=&
(\delta([{}^{\beta_k} \alpha_{r+1},{}^{\beta_k}\alpha_1],
\ldots,\delta([{}^{\beta_k} \alpha_{r+1},{}^{\beta_k}\alpha_r]))\nonumber \\
&=&
(\delta([ \alpha_{r+1},\alpha_1],
\ldots,\delta([ \alpha_{r+1},\alpha_r]))\circ B_k\\
&=& \beta_k^{-1}(\tau(\delta)),\nonumber\end{eqnarray}
for all $\delta \in C^1(\Pi,V_\lambda^\vee)$ and $k=1,\ldots,r;$
where the first and the last equality hold by definition.
Equality  (\ref{3.2})
 follows from an elementary but  tedious
 computation, using the crossed homomorphism 
relation and  the action of $\beta_k$ on $(\alpha_1,\ldots,\alpha_r).$ 
\Endproof

\section{Convolution of local systems}\label{s2}

It is the aim of this section to give an interpretation 
of the multiplicative
version of the convolution in terms of the cohomology 
of local systems
on the punctured sphere.\\

\subsection{Local systems}\label{localsyst}

Let $W$ be a connected topological 
manifold. 
A  {\em (complex) local system of rank $n$} 
on $W$   is a 
sheaf $\FF$ of complex vector spaces which is locally 
isomorphic to the constant sheaf $\CC^n.$ 
The category of local systems on $W$ is denoted by 
$\LocSys(W).$ It is closed under 
tensor product and taking duals.
The dual local system  
of a local system  $\FF$ on $W$ is denoted by $\FF^\vee.$ 
The stalk of $\FF$ at $w_0\in W$ is denoted by $\FF_{w_0}.$\\
If $\gamma$ is a closed path in $W$ starting at $w_0,$ then there exists 
a unique linear transformation $\Mon(\gamma)$ such that 
the stalk $ \gamma^*(\FF)_1$ is canonically 
isomorphic to $ \Mon(\gamma)\cdot  \gamma^*(\FF)_0.$ 
Composition
 of paths gives rise to the monodromy 
representation (see \cite{Deligne70}): 
$$ \Mon: =\Mon(\FF):\pi_1(W,w_0)\to \GL(\FF_{w_0}).$$ 
It is well known, that the functor
$$ \LocSys(W)\to \Mod(\CC[\pi_1(W,w_o)]),\,\, \FF \mapsto \FF_{w_o}$$
is an equivalence of categories,
see \cite{Deligne70}, Corollaire 1.4.

\subsection{The middle convolution functor $MC_\lambda$
for  local systems}\label{42}

\begin{defn}\label{Kummer}{\rm  Let $\gamma$ be a closed path in 
$\CC^\times$ which has  initial point $x_0$ and 
encircles $0$ once in counterclockwise direction,
$\lambda \in \CC^\times$ and 
$$\chi: \pi_1(\CC^\times
,x_0)\to \GL(\CC),\, \gamma \mapsto \lambda.$$
The local system on $\CC^\times$
corresponding  to the (module associated to the)
homomorphism $\chi$ is  called the 
{\em Kummer sheaf} associated to $\lambda$ and $y_0$ and is 
denoted by 
 $\L_\lambda.$  }
\end{defn}

Let $\alpha_1,\ldots,\alpha_{r+1}$ (resp. $\beta_1,\ldots,\beta_r$) 
be as in the previous subsections and
  $\F$  be the 
local system associated to the representation
$$ \rho:\pi_1(X,x_0)\to \GL(V),\,\, \alpha_i\mapsto A_i,\,i=1,\ldots , r.$$
 For 
 $\lambda\in \CC,$ let 
 $C_\lambda(\A)=(B_1,\ldots,B_r)$ $\in \GL(V^r)^r$
and $MC_\lambda(\A)=(\tilde{B}_1,\ldots,\tilde{B}_r)\in 
\GL(V^r/(\K+\L))^r.$
We define $C_\lambda(\FF)$ to be the 
local system associated to 
$$ \pi_1(X,x_0)\to \GL(V^r),\,\, \beta_i\mapsto B_i,\,i=1,\ldots , r.$$
Similarly, let $MC_\lambda(\FF)$ be the local system associated 
to 
$$ \pi_1(X,x_0)\to \GL(V^r/(\K+\L)),\,\, \beta_i\mapsto \tilde{B}_i,\,i=1,\ldots , r.$$

\begin{prop} The local system $C_\lambda(\FF),$
resp. $MC_\lambda(\FF),$ is canonically isomorphic
to the local system corresponding to the representation
$$ \pi_1(X,x_0)\to \GL(V),\,\, \alpha_i\mapsto B_i,\,i=1,\ldots , r,$$
resp.
$$ \pi_1(X,x_0)\to \GL(V),\,\, \alpha_i\mapsto \tilde{B}_i,
\,i=1,\ldots , r.$$
\end{prop}

\proof There exists, up to homotopy, a
unique path $\gamma$ in $X$  with initial point 
$y_0$ and endpoint $x_0$ such that $\beta_i=\gamma^{-1}\alpha_i\gamma$
for $i=1,\ldots,r.$
This path induces canonical isomorphisms 
$$C_\lambda(\FF)|_{y_0}\to C_\lambda(\FF)|_{x_0}, $$
resp. 
$$MC_\lambda(\FF)|_{y_0}\to MC_\lambda(\FF)|_{x_0}, $$
which are compatible with the induced  isomorphism of 
fundamental groups
$$ \pi_1(X,y_0)\to \pi_1(X,x_0),\,\, \beta_i \mapsto \alpha_i=
\gamma\alpha_i\gamma^{-1}$$ and the action of the two fundamental groups
on their respective fibers.
\Endproof

In view of Proposition \ref{deli2} one obtains covariant, end-exact, functors
$$C_\lambda: \LocSys(X)\to \LocSys(X), \,\,\FF \mapsto C_\lambda(\FF)$$
and 
$$MC_\lambda: \LocSys(X)\to \LocSys(X), \,\, \FF \mapsto MC_\lambda(\FF).$$
Moreover, all the properties of $MC_\lambda,$ as given in Theorem 
\ref{eigen}, immediately translate into the language of local systems.\\

The following definition is justified by the results of the
next 
subsection:

\begin{defn}{\rm  
The local system 
$C_\lambda(\FF)$ (resp. 
$MC_\lambda(\FF)$) is called 
the {\em convolution} (resp. {\em middle convolution})
of $\FF$ with $\L_\lambda.$}
\end{defn}

\subsection{Cohomological interpretation of $MC_\lambda$}\label{cc}

Let  
$\S$ a sheaf of complex vector spaces on $W.$ An $i$-cochain $\psi$
is a map which  associates to
any $i$-chain $\sigma:\Delta^i\to W$  an element $\psi(\sigma)\in \CC.$
The set of $i$-cochains is denoted by $C^i(W).$ 
Consider the (injective and torsionless, see \cite{Warner71})
resolution of the constant sheaf $\CC$ on $W$
via cochains
$$0\to\CC\to C^0(W)\to C^1(W)\to C^2(W)\to \ldots$$
and let $H^i(W,\S):=H^i(\Gamma(C^*\otimes \S)).$ \\

Let $p: W_1\to W_2$ be a continuous map
of topological spaces
and $\S$ a sheaf  on $W_1.$
The sheaf associated to the presheaf
$$ U \mapsto S_U:=H^i(p^{-1}(U), \S|_U)\quad 
\mbox{\rm ($U$ open in $W_2$)}$$
is denoted by $R^i p_*(\S)$
(it is well known that 
$R^i p_*$ can be viewed as  the  $i$-th higher direct image functor 
 of $p_*$).\\

Let
$$E=\{(x,y)\in \CC^2 \mid x,y\not= t_i,\, i=1,\ldots,r,\, x\not= y\},$$
 ${\rm p}_i:E\to X,\, i=1,2,$ be the $i$-th projection, $$q:E\to \CC^\times,\,
(x,y)\mapsto y-x,$$ $j: E \to \PP^1(\CC)\times X$ the tautological 
inclusion
 and 
$\bar{p}_2: \PP^1(\CC)\times X \to X$ the  (second) projection onto $X.$ \\

\begin{thm}\label{ccc} Let $\F$ be a 
local system on $X,$ $\lambda \in \CC^\times\setminus 1$
and $\L_\lambda$ the Kummer 
sheaf associated to $\lambda.$ Then 
$$  MC_\lambda(\F)\cong  R^1 (\bar{p}_2)_*(j_*({\rm p}_1^*(\F)\otimes 
q^*(\L_\lambda))).$$
 \end{thm}

\proof 
Let $\F$ be the local system associated to 
a representation
$$ \pi_1(X,x_0)\to \GL(V),\,\, \alpha_i \mapsto A_i,\, i=1,\ldots,r,$$
where $\alpha_i$ is as in Section \ref{localsyst}.
Let  $\G:={\rm p}_1^*(\F)\otimes 
q^*(\L_\lambda)$ and  
$\G_{y_0}$  the restriction of $\G$ to $X(y_0)$
(thus $\G_{y_0}$ corresponds to the $\Pi$-module 
$V_\lambda$ of the last section).
Let $\psi\in C^1(X(y_0),\G_{y_0}^\vee)$ be a closed cochain
 and $\sigma_1,\,\sigma_2$ closed paths in 
$X(y_0),$ based at $x_0.$
 By definition,
 $$ \psi(\sigma_1\sigma_2)|_0=\psi(\sigma_2)|_0
+\Mon(\sigma_2)^{-1}\psi(\sigma_1)|_0.$$ 
This induces an isomorphism
$$ H:H^1(X(y_0),\G_{y_0}^\vee)
\to H^1(\Pi,V_\lambda^\vee),\,\, [\psi] \mapsto 
[(\sigma \mapsto \psi(\sigma)|_{0})].$$ 

Since ${\rm p}_2:E\to X$ is a locally trivial fibration,
$R^1({\rm p}_2)_*(\G^\vee)$ is a local system. By construction, 
 the monodromy action of $\beta_k\in 
\pi_1(X,y_0)$ on $R^1({\rm p}_2)_*(\G^\vee)|_{y_0}=
H^1(X(y_0),\G_{y_0})$
is the one which is induced 
by sending $\alpha_i$ to ${}^{\beta_k^{-1}}\alpha_i.$
 This yields a canonical  isomorphism of
 $\pi_1(X,y_0)$-modules between $H^1(X(y_0),\G_{y_0}^\vee)$
and $H^1(\Pi,$ $V_\lambda^\vee).$
Thus, by
 Theorem \ref{Bk},
one has 
a canonical isomorphism
\begin{eqnarray}  R^1({\rm p}_2)_*(\G^\vee) \cong C_\lambda(\F)^\vee.\end{eqnarray}

 Consider the following subspaces of
 $(V_\lambda^\vee)^r= C_\lambda(\F)^\vee|_{y_0}:$
 $$V_1:=(\im((A^\vee_1)^{-1}-1),\ldots,\im((A^\vee_r)^{-1}-1))$$
and
\begin{eqnarray} V_2:&=&\{(w_1,
 w_2 ,\ldots,
w_r) \in (V_\lambda^\vee)^r \,\mid \,\nonumber\\
&& \quad\quad\quad
(\sum_{i=1}^{r-1}((A_{i+1}\cdots A_r)^\vee)^{-1}w_i)+w_r 
\in \im((A^\vee_1\cdots A^\vee_r )^{-1}\lambda -1)\}. \nonumber 
\end{eqnarray}
One can easily check that $V_1\cap V_2\leq (V_\lambda^\vee)^r$ is the 
$\pi_1(X,y_0)$-submodule  which corresponds 
to $MC_\lambda(\F)^\vee.$ 

The image 
of the cohomology with compact 
supports $H^1_c(X(y_0),\G^\vee_{y_0})$
in $ H^1(X(y_0),$
$\G^\vee_{y_0})$ 
is mapped under $ \tau\circ H$ isomorphically 
onto $V_1\cap V_2.$
This can be seen using similar 
argu\-ments as Shimura \cite{Shimura71}, Chap. 8,
or by writing $\alpha_i$ as the product 
$\bar{\gamma}_i^{-1}\hat{\gamma}_i\bar{\gamma}_i$ (where $\bar{\gamma}_i$
is a  path which starts at $x_0$ and goes near to the singularity 
$t_i,$ and $\hat{\gamma}_i$ moves along a small circle around $t_i$)
and using the compact supports condition at $t_1,\ldots,t_k$ and $\infty.$
The image of 
$H^1_c(X(y_0),\G^\vee_{y_0})$ in 
$ H^1(X(y_0),\G^\vee_{y_0})$ is canonically isomorphic to 
$H^1(j_*(\G^\vee_{y_0}))$ (see \cite{Looijenga97}, Lemma 5.3). 
Therefore,
$$ MC_\lambda(\F)^\vee \cong R^1(\bar{p}_2)_*(j_*(\G^\vee)).$$

Finally, the Poincar\'e pairing yields an isomorphism 
$$R^1(\bar{p}_2)_*(j_*(\G^\vee))^\vee \cong R^1 (\bar{p}_2)_*(j_*(\G))= 
R^1 (\bar{p}_2)_*(j_*({\rm p}_1^*(\F)\otimes 
q^*(\L_\lambda)))$$ (see e.g. \cite{Looijenga97}, Lemma 5.3).
\Endproof

\begin{rem}\label{bas} {\rm  i) The resolution via singular cochains allows one to use 
ground fields  different 
from $\CC$ as coefficients of 
cohomology. One could even work in the category of 
local systems over
principal ideal domains, see \cite{Warner71}.}
\end{rem}

\section{The middle convolution transformation $mc_\mu$ 
of Fuchsian systems}

\subsection{Definition of $mc_\mu$ for tuples of matrices}\label{addi}

In this section we recall the additive convolution as given in 
\cite{DR}, App. A.\\

Let $K$ be any field and 
 $\a=(a_1,\ldots,a_r),\, a_k \in K^{n \times n}.$
 For $\mu \in K$ one can define
blockmatrices $b_k,\, k=1,\ldots,r,$ as follows: 
  \[ b_k:= \left( \begin{array}{ccccccc}
                   0 &  &  & \ldots& & 0\\
                   & \ddots &  & & &\\
               a_1 & \ldots\,\, a_{k-1}&  a_k+ \mu  &  a_{k+1} & \ldots 
&  a_r \\
 %            0 & \ldots & & 0 & 0 &   \ldots  \\     
               &   &  &  \ddots & &   \\             
                   0 &    & \ldots& &  & 0
          \end{array} \right) \in K^{nr \times nr},\]
where $b_k$ is  zero outside the $k$-th block row.\\

There are  the following left-$\langle b_1,\ldots,b_r \rangle$-invariant 
subspaces of the 
column vector space $K^{nr}$ (with the tautological action of 
$\langle b_1,\ldots,b_r \rangle$):\\

\[ {\frak{k}}_k = \left( \begin{array}{c}
          0 \\
          \vdots \\
            0\\
          \ker(a_k) \\
           0 \\
           \vdots \\
           0 
        \end{array} \right) \quad \mbox{($k$-th entry)},\, k=1,\dots,r,\]
and 
\[     {\frak{l}}=\cap_{k=1}^r {\rm ker}(b_k)
={\rm ker}(b_1+\ldots+b_r).
\]
Let ${\frak{k}}:=\oplus_{k=1}^r {\frak{k}}_k.$

If $\mu \not=0$ then 
$${\frak{l}}= \langle\left( \begin{array}{c}
                       v \\
               \vdots \\
                       v
            \end{array} \right)\mid  v \in \ker(a_1+ \cdots+a_r+\mu) \rangle.$$
and 
$${\frak{k}}+{\frak{l}}={\frak{k}}\oplus {\frak{l}}.$$
We fix an isomorphism $I$ 
between $K^{nr}/({\frak{k}}+{\frak{l}})$ and $K^m.$

\begin{defn} {\rm We call $c_{\mu}(\a):=(b_1,\dots,b_r)$ the {\em (additive
version of the)
 convolution}
of $\a=(a_1,\ldots,a_r)$ with $\mu.$
The tuple of matrices 
$mc_{\mu}:=(\tilde{b}_1,\dots,\tilde{b}_r)\in K^{m\times m}, $
where $\tilde{b}_i$ is induced by the action of $b_i$ on 
$K^m(\simeq K^{nr}/({\frak{k}}+{\frak{l}})),$ is called the
{\em (additive version of the) middle convolution} of $\a$ with $\mu.$}
\end{defn}

\subsection{The definition of 
$mc_\mu$ for Fuchsian systems and monodromy of 
differential systems}\label{fuchs}

Let $T:=\{t_1,\ldots,t_r\},\,X:=\CC\setminus T$ 
and ${\bf c}:=(c_1,\ldots,c_r),\, c_i \in \CC^{k\times k}.$
The Fuchsian system 
$$Y'=\sum_{i=1}^r \frac{c_i}{x-t_i} Y$$ is denoted by $D_{\bf c}.$

\begin{defn}{\rm  Let $\a:=(a_1,\ldots,a_r),\, a_i \in \CC^{n\times n},$
and $\mu \in \CC.$ The Fuchsian system $D_{c_\mu(\a)}$ (resp.
$D_{mc_\mu(\a)}$) is called the {\em convolution (resp. middle 
convolution)} of $D_\a$ with $\mu.$}
\end{defn}

Let $\gamma_1,\ldots, \gamma_{r+1}$ be a homotopy base
of $\pi_1(X,o),$ $D$ a linear system of differential equations
which has no singularities in $X$ and $F$ a fundamental system of $D,$
consisting of functions which are defined in a small neighborhood 
of $o.$  
Analytic continuation of $F$ along $\gamma_i$ transforms 
$F$ into  $F\cdot \Mon(\gamma_i).$
We call the tuple $$\Mon(D):=(\Mon(\gamma_1),\ldots,\Mon(\gamma_r))$$ the 
tuple of {\it monodromy generators} of $D$ with respect to 
$F$ and $\gamma_1,\ldots,\gamma_r.$ 

\begin{rem} \label{rightleft}{\rm i) An element $\gamma_i \in  \pi_1(X,o)$
acts (via $\Mon(\gamma_i)$) 
from the right on the vector space ${\cal S}$ spanned 
by the rows of the fundamental system $F$ of $D.$ 
Let $\F$ denote the local system $\F$ formed by the  solutions of $D$
(locally at $o$ given by the columns of $F$) and fix the isomorphism 
$$\F_o\to \CC^n, f_i(o) \mapsto e_k,$$ where $f_i$ denotes the 
$i$-th column of $F$ and 
$e_k$ is the $k$-th standard vector of $\CC^n.$
Then  the monodromy of $\F$ with respect to 
$\gamma$ is given 
by the same matrix $\Mon(\gamma_i)$
acting from the left on $\F_0\simeq \CC^n.$

ii) As a factor system of $D_{c_\mu(\a)},$ the middle convolution 
$D_{mc_\mu(\a)},\, mc_\mu(\a)\in (\CC^{m\times m})^r,$ 
can be constructed by a base change, transforming 
a basis of ${\frak k}+{\frak l}$
to the first $nr-m$ standard vectors, and cutting out 
the $m\times m$-block matrices corresponding to the last 
$m$ entries. The same construction applies for a fundamental
matrix of $D_{c_\mu(\a)}.$ We say that a matrix 
whose columns
are solutions of  $D_{c_\mu(\a)}$ (not necessarily a fundamental matrix
of $D_{c_\mu(\a)}$)
{\em gives rise to a fundamental matrix of} 
$D_{mc_\mu(\a)},$ if the resulting matrix under the above construction
of transforming and cutting out is a fundamental matrix of
$D_{mc_\mu(\a)}.$
}\end{rem}

\section{Compatibility of $MC_\lambda$ and $mc_\mu$}\label{Eule}

In this section we relate the additive version of the 
convolution to the multiplicative version (Subsection \ref{period}).

\subsection{The Euler transform}\label{cont}

A commutator 
$$[\alpha_{i},\alpha_j]=
\alpha_{i}^{-1} \alpha_{j}^{-1} \alpha_{i}  \alpha_{j} $$
is called a {\em Pochhammer contour}.
Pochhammer contours are widely used in the theory of 
ordinary differential equations, see \cite{Ince56}, \cite{Kohno99}
 and \cite{Yoshida97}.

\begin{defn}{\rm   Let $ \mu \in \CC,$ 
  $g:=(g_{i,j})$ be a matrix whose entries 
$g_{i,j}$ are (multi-valued) 
functions which are   holomorphic 
 on $X.$ 
The path $\alpha_{r+1}$ encircles an open neighbourhood $U$
of $y_0.$
 The matrix valued function
 \[   I_{[\alpha_{r+1},\alpha_i]}^\mu ( g)(y) := \int_{[\alpha_{r+1},\alpha_i]} g(x) 
(y-x)^{{\mu-1}}   dx ,\,\, y \in U,\]
is called the {\em Euler transform} of $g$ with respect to 
$[\alpha_{r+1},\alpha_i]$ 
and 
$\mu$.} \end{defn}

The next lemma shows that the Euler transformation
is compatible with the convolution:

\begin{lemma}\label{Poch} Let 
$\a:=(a_1,\ldots,a_r),\, a_i \in \CC^{n\times n},$
and $\mu_1,\, \mu_2 \in \CC.$ 
 If  $g(x)$ is a solution of $D_{c_{\mu_1}(\a)},$ 
 then $I_{[\alpha_{r+1},\alpha_i]}^{\mu_2} (g)(y)$
 is a solution for 
 $D_{c_{\mu_1+\mu_2}(\a)},$ where $y$ is contained 
in an open neighborhood of $y_0$ which is encircled 
by $\alpha_{r+1}.$
\end{lemma}

\proof In the following, we  omit 
the subscript $[\alpha_{r+1},\alpha_i]$ at the integral sign.
For $y \in U,$ 
 \begin{eqnarray}
(y-T)  \frac{d(I_{[\alpha_{r+1},\alpha_i]}^{\mu_2} ( g))}{dy}& = &
(y-T)
\int\frac{d}{dy} g(x)(y-x)^{\mu_2-1}dx\nonumber\\
&=  & \int ((y-x) + (x-T)) (\frac{d}{dy}g(x)(y-x)^{\mu_2-1})dx\nonumber\\
&=& (\mu_2-1)I_{[\alpha_{r+1},\alpha_i]}^{\mu_2}(g)(y)\\&&+\,(\mu_2-1)\int(x-T)g(x)
(y-x)^{\mu_2-2}dx \label{2.1},
\nonumber \end{eqnarray}
where one is allowed to 
differentiate under the integration sign since $[\alpha_{r+1},\alpha_i]$
is compact. 
One has 
 \begin{eqnarray}
0&=& \int \frac{d}{dx}((x-T)g(x)(y-x)^{\mu_2-1})dx \label{conto}\\
&=& \int g(x)(y-x)^{\mu_2-1}dx +\nonumber\\&&
\int (x-T)g'(x)(y-x)^{\mu_2-1}dx -\nonumber\\&&
(\mu_2-1)\int(x-T)g(x)(y-x)^{\mu_2-2}dx.\nonumber\end{eqnarray}
Therefore
 \begin{eqnarray}(\mu_2-1)\int(x-T)g(x)(y-x)^{\mu_2-2}dx& =& 
\int g(x)(y-x)^{\mu_2-1}dx +\nonumber\nonumber  \\&&
\int (x-T)g'(x)(y-x)^{\mu_2-1}dx.\nonumber\end{eqnarray}
Using the last equality one sees that 
\begin{eqnarray}(5)&
=& (\mu_2-1)I_{[\alpha_{r+1},\alpha_i]}^{\mu_2}(g)+I_{[\alpha_{r+1},\alpha_i]}^{\mu_2}(g)+
\int(x-T)g'(x)(y-x)^{\mu_2-1}dx\nonumber \\
&
=&
(\sum_{k=1}^r b_k)I^{\mu_2}_{[\alpha_{r+1},\alpha_i]} (g)(y)\nonumber, \end{eqnarray}
where $c_{\mu_1+\mu_2}(\a)=(b_1,\ldots,b_r)$
(use that 
 $g$ is a solution of $c_{\mu_1}(D_\a)$).
\Endproof

\begin{rem}\label{contrast}{\rm  
The use of Pochhammer contours 
is an essential ingredience in the proof of the above lemma
(see Formula (\ref{conto})).}
\end{rem}

In the following, $F$ denotes a fundamental system of a 
Fuchsian system $D_\a$ and 
$$G(x):=\left(\begin{array}{c}F(x)(x-t_1)^{-1}\\
\vdots\\
F(x)(x-t_r)^{-1}\end{array}\right).$$

The next results will be used in the proof of Theorem \ref{monodromy}:

\begin{lemma} \label{Pochloop}
i) The columns of $G$ are solutions of
 $D_{c_{-1}(\a)}.$

ii)
 $I_{[\alpha_{r+1},\alpha_i]}^{\mu}(G)=  I_{\alpha_{i}}^{\mu}(G)(1-e^{2 \pi i \mu})
-I_{\alpha_{r+1}}^{\mu}(G) (1-\Mon(\alpha_i)).$
\end{lemma}

\proof
 The first assertion follows from a straightforward computation.
The second assertion follows from the definition of $[\alpha_{r+1},\alpha_i],$
using
 the effect of the monodromy on the integrand, see \cite{Ince56}, Chap.
18.
\Endproof

\begin{cor}\label{Cauchy}
i) If $\mu$ is a positive integer, then
 \[I_{[\alpha_{r+1},\alpha_i]}^{\mu}(G)=0. \]

ii) If $\mu=0$ or a negative integer, then
  \[I_{[\alpha_{r+1},\alpha_i]}^{\mu}(G)= {2 \pi i \over -\mu! } G^{(-\mu)}(y) (-1+\Mon(\alpha_i)). \]
\end{cor}

\proof 
The claims follow from the above lemma and Cauchy's integral formula. 
\Endproof

\begin{lemma} \label{add}
 Let 
 $ Y'= \sum { a_i \over x-t_i} Y $ be a Fuchsian system
 with a nontrivial monodromy group
 and  $\mu \not \in \ZZ.$
 Then there exits an $i$ and a solution $f(x)$ such that
 \[ \int_{[\alpha_{r+1},\alpha_i]} {f(x) \over x-t_i} (y-x)^{\mu-1}dx \neq 0. \]
\end{lemma}

\proof
 We can assume that we have nontrivial
 monodromy at $t_1=0.$
 If the monodromy is not unipotent then we can find an
 entry 
 \[g(x)=x^\alpha \sum_{j=0}^\infty x^j a(j) ,\,\alpha \not\in \ZZ,
\, a(0)\not=0 \]
 of a  solution $f(x)$ near $t_1.$
 Then
$$I_{[\alpha_{r+1},\alpha_1]}^\mu(g)=
\sum_{j=0}^\infty a(j) \int_{[\alpha_{r+1},\alpha_1]}x^{\alpha+j}
(y-x)^{\mu-1} dx .$$
Using similar arguments as in \cite{Yoshida97}, Chap. IV,
one can prove that 
$$\int_{[\alpha_{r+1},\alpha_1]}x^{\alpha+j}
(y-x)^{\mu-1} dx =y^{\alpha+j+\mu}\beta(\alpha +i,\mu),$$
where $\beta(\alpha +i,\mu)\not=0,$ if $\alpha +i,\mu\not\in  \ZZ.$

 In the case of nontrivial unipotent monodromy
 at $t_1=0,$
 there exists a solution
which has an entry $g(x)=h_0(x)+\log(x)h_1(x)$ near $0,$ where
 $h_0,h_1$ are meromorphic at $0$ and $h_1 \neq 0$ (see \cite{Ince56}, 16.3).
The claim follows now from an easy exercise, using 
integration
by parts, Lemma \ref{Pochloop} and Corollary \ref{Cauchy}.
\Endproof

\subsection{The Riemann-Hilbert correspondence for 
 $MC_\lambda$}\label{period}

In the notation of the preceding sections.
Let  $\a:=(a_1,\ldots,a_r),\, a_i \in \CC^{n\times n}.$ Let 
  $F$  be a fundamental matrix  of  the Fuchsian system 
$D_{\a}:Y'=\sum \frac{a_i}{x-t_i}Y$ and  
         $$G(x):=\left(\begin{array}{c}F(x)(x-t_1)^{-1}\\
\vdots\\
F(x)(x-t_r)^{-1}\end{array}\right).$$

\begin{defn}{\rm  Let 
$\mu \in \CC.$ The matrix 
$$I^{\mu}:=I^\mu(y):=
(I^\mu_{[\alpha_{r+1},\alpha_1]}(G)(y),\ldots, I^\mu_{[\alpha_{r+1},
\alpha_r]}(G)(y)) $$ is called 
the {\em period matrix}.} \end{defn}

\noindent {\bf Remark:} It follows from the definitions  
 that if the period 
matrix $I^\mu$ is invertible, then 
it  describes 
the pairing between $H_1(X(y_0),\G_{y_0})$ and 
$H^1_{DR}(X(y_0),\G_{y_0}^\vee),$ where 
$X(y_0)$ and $\G$ are as in the proof of 
Theorem \ref{ccc}.
 In the next theorem we will give 
criteria for $I^\mu$ to be invertible, i.e., the rows of
$I^\mu(y_0)$ exhibit a base of  
$H^1_{DR}(X(y_0),\G_{y_0}^\vee).$

In a similar way as described in 
Yoshida \cite{Yoshida97}, Chap. iv, it can be shown that the matrix 
${\frak H}$ which occurs in Lemma \ref{form2} has a natural interpretation
as an intersection matrix of (``loaded'') cycles 
$c\in H_1(X(y_0),\G_{y_0}).$\\

The next theorem shows the relation between the additive 
and multiplicative versions of the convolution:

\begin{thm}\label{monodromy} Let 
$\a:=(a_1,\ldots, a_r),\, a_i\in \CC^{n\times n},$
 $\Mon(D_\a)
=(A_1,\ldots,A_r)\in \GL_n(\CC)^r$ its tuple of 
monodromy generators,  $\mu \in \CC\setminus \ZZ$
and $\lambda:=e^{2\pi i \mu}.$ If 
the generated subgroup $\langle A_1,\ldots,A_r \rangle$ is an 
irreducible subgroup of $\GL_n(\CC)$ 
and if at least two elements of $A_1,\ldots,A_r$ are $\not= 1,$ then
the following statements hold:\\

i) The columns of the period 
matrix $I^{\mu}(y)$ are solutions of $D_{c_{\mu-1}(\a)},$
where $y$ is contained in a small open neighbourhood $U$ of $y_0.$\\

ii) For $v_i \in \ker(A_i-1),\, i=1,\ldots,r,$ 
(resp. $v \in \ker(A_1\cdots A_r\lambda
-1)$) assume that  the residues
of $G(x)v_i$ at $t_i$  (resp. the residues of 
$x^{\mu-1}G(x)v$ at  $\infty$)
is not identically zero. Then the 
period matrix $I^{\mu}(y),\, y\in U,$ is a 
fundamental matrix  of $D_{c_{\mu-1}(\a)}.$
Further,
the tuple of
monodromy generators of $D_{c_{{\mu-1}}({\a})}$ with respect
to $I^{\mu}(y)$ and the paths $\beta_1, \ldots, \beta_r$
is       $ C_{\lambda}(\Mon(D_{\a})), $ i.e.,
     \[   \Mon(D_{c_{{\mu-1}}(\a)})=C_{\lambda}(\Mon(D_\a)). \]

iii) Assume that  
$$\rk(a_i) = \rk(A_i-1)\quad \mbox{\rm  and}\quad  
\rk(a_1+\cdots+ a_r+\mu)=\rk( \lambda\cdot A_1\cdots A_r-1).$$
The matrix $I^{\mu}(y)$ gives rise to 
 a fundamental matrix $\tilde{I^\mu}(y),\, 
y\in U,$
of the system  
$D_{mc_{{\mu-1}}(\a)}$ (see Remark \ref{rightleft}, ii)).
The tuple of monodromy 
generators of $D_{mc_{{\mu-1}}({\a})}$ with respect to
  $\tilde{I^{\mu}}(y)$ and the paths $\beta_1, \ldots, \beta_r$
is       $ MC_{\lambda}(\Mon(D_{\a})), $ i.e.,
     \[   \Mon(D_{mc_{{\mu-1}}(\a)})=MC_{\lambda}(\Mon(D_\a)). \]
\end{thm}

\begin{rem}{\rm 
a)  It follows from the proof
that one can weaken the  assumptions of Theorem \ref{monodromy}
such that tuple $(A_1,\ldots,A_r)$ 
fulfills 
the conditions $(*)$ and $(**)$ of Subsection \ref{property}
instead of the irreducibility and non-triviality condition on 
$A_1,\ldots,A_r.$

b) In Theorem \ref{monodromy}, iii), if $\rk(a_i) > \rk(A_i-1)$
then the differential system, which corresponds to 
(the local system corresponding to) $MC_{\lambda}(A_1,\ldots, A_r)$
is a factor system of $D_{mc_{{\mu-1}}(\a)}.$
 }\end{rem}

\noindent {\bf Proof} of i): This follows from Lemma \ref{Poch} and 
Lemma \ref{Pochloop} i).\Endproof

\noindent {\bf Proof}  of ii): Let 
$V_\lambda$ denote the $\pi_1(X\setminus
\{y_0\},x_0)$-module whose underlying vector space is 
the column vector space $\CC^n$ on which  
$\alpha_i$ acts via $A_i$ ($1\leq i \leq r$) and 
$\alpha_{r+1}$ acts via $\lambda.$ Let $\CC_n$ denote the 
space of row vectors and $V_\lambda^\vee$ the dual module 
of $V_\lambda$ with respect to 
$\CC_n\times \CC^n \to \CC, \,\, (w,v)\mapsto w\cdot v.$
Let $F^j$ denote the $j$-th row
of $F$ and 
$$ \delta_{i,j}: \pi_1(X\setminus
\{y_0\},x_0)\to V_\lambda^\vee,\, \gamma \mapsto 
 \int_{\gamma}{F^j(y_0-x)^{\mu-1}}\frac{dx}{x-t_i}.$$ 
By the 
properties of the integral,
  $\delta_{i,j}$ is an 
element in $C^1(\pi_1(X\setminus
\{y_0\},x_0),V_\lambda^\vee).$
It follows that the rows of $I^\mu(y_0)$ are exactly the images 
of the 
crossed homomorphisms $\delta_{i,j}$ under the twisted evaluation map.
By the definition of $I_{[\alpha_{r+1},\alpha_i]}^\mu,$ analytic 
continuation of $I^\mu(y)$ along the 
path $\beta_k$ transforms $I_{[\alpha_{r+1},\alpha_i]}^\mu(G)$
into $I_{[{}^{\beta_k}\alpha_{r+1},{}^{\beta_k}\alpha_i]}^\mu(G).$
 It follows then from Theorem \ref{Bk}
(Formula (\ref{3.2})), that
the matrix which 
 describes this  transformation is the matrix $B_k,$ where 
$$C_\lambda(A_1,\ldots,A_r)=(B_1,\ldots,B_r).$$ 

In order to prove ii),
it remains to prove that the columns of $I^\mu(y)$ form a 
fundamental set of solutions. This follows from the Lemmata below:\\

Consider the vector space of solutions
$J:=I^{\mu}(y)\cdot \CC^{nr},$
with $y$ in a small neighborhood 
of $y_0.$  Let further 
 $\K_i,\,\K$ and $ \L$ be as
in Subsection \ref{katz} and   $\hat{\K}_i:=I^\mu(y)\cdot \K_i,\,
\hat{\K}:=I^\mu(y)\cdot\K$ and $\hat{\L}:=I^\mu(y)\cdot\L.$

\begin{lemma}\label{kern} The kernel of the map
$$I^\mu:\CC^{nr} \to J,\, (v_1,\ldots,v_r)^\tr\to I^\mu(y)\cdot
(v_1,\ldots,v_r)^\tr $$ is a $\langle B_k:
k=1,\ldots,r\rangle$-module. 
\end{lemma}
\proof If $I^{\mu}v=0,$ then $I^{\mu}B_kv=0.$
\Endproof

If $G=(g_{i,j}(t))$ is a vector valued function which is 
componentwise meromorphic a $t_k,$ then
$\Res_{t_k}(G)$ denotes the vector of 
residues $(\Res_{t_k}(g_{i,j}(t))).$

\begin{lemma}\label{Resi}
Let $ \mu \not \in \ZZ.$ Then 
 the functions in $\hat{\K}_i$ (resp. $\hat{\L}$) have at most a
 singularity at $t_i$ (resp. $\infty$). Moreover,

i) \[  \hat{\K}_i= \langle \Res_{t_i}(G(x)v) (y-t_i)^{\mu-1}\mid
                v\in \ker(A_i-1) \rangle \]
 for $i=1,\ldots,r.$
 
ii)  \[ \hat{\L}= \langle \Res_{\infty}(x^{\mu-1}G(x)v)\mid 
                    v\in \ker(A_1\cdots A_r \lambda -1) \rangle \]           
\end{lemma}

\proof
 i) One has
    \begin{eqnarray*}
      \hat{\K}_i =  I^{\mu} \K_i & = & I^{\mu}_{[\alpha_{r+1},\alpha_i]}(G(x)) \ker(A_i-1) \\
     &=& I^{\mu}_{\alpha_i}(G(x)) \ker(A_i-1)
     \end{eqnarray*}
    by Lemma~\ref{Pochloop}.
    The claim follows from  Cauchy's integral formula since
$G(x)v$ (as  matrix valued function in $x$)
    is meromorphic at $t_i$ for $v \in \ker(A_i-1).$

ii) Using Lemma~\ref{Pochloop} one easily sees that
    $$\hat{\L} =I^{\mu}_{[\alpha_{r+1},\alpha_{\infty}]}(G(x)) 
\ker(A_1\cdots A_r \lambda-1),$$
    where $\alpha_{\infty}= \alpha_1 \cdots \alpha_{r+1}.$ 
    Using the same arguments as in i) the claim follows.
\Endproof

\begin{cor}\label{hut}
 One has
 \[ \hat{\K}+\hat{\L}= \oplus_i \hat{\K}_i \oplus\hat{\L} \] as
 a left-$\langle B_1,\ldots,B_r \rangle$-module.
\end{cor}

\begin{lemma}\label{4.14}
 If the conditions $(*)$ and $(**)$ of Subsection 
\ref{property} hold for $\Mon(D_\a)=(A_1,\ldots,A_r),$ then
 \[ \ker(I^{\mu}) \leq {\K}+{\L}. \]
\end{lemma}

\proof
Assume that  $\ker(I^{\mu}) \not\leq \K+\L.$
Let $O \leq V_1 \leq \ldots \leq V_k=V_\lambda$ be a composition series of 
$V_\lambda$
(as a module). Let further $V_i^{r}$ be the corresponding (diagonal)
subspace of 
$V_\lambda^r=\CC^{nr}$ and
$\tilde{V}_i:= V_i^r+\K+\L \mod \K+\L.$
It follows from Theorem \ref{eigen}, iii), and 
\cite{DR}, Lemma 2.8,  that
 $O \leq \tilde{V}_1 \leq \ldots \leq \tilde{V}_k=V_\lambda^r/(\K+\L)$
is a composition series of $V_\lambda^r/(\K+\L)$ 
(as $\langle B_1, \ldots, B_r \rangle$-module).
Since $\ker(I^{\mu})$ is a $\langle B_1, \ldots, B_r \rangle$-module,
there exists a  
 $\langle A_1,\ldots,A_r \rangle$-module $W \leq V_\lambda$ such that
$W^r+\K+\L\leq \ker(I^{\mu}) +\K+\L.$
We assume  that $W$ is  minimal and nontrivial. 
By minimality, $(*)$ and $(**)$ also hold 
for $W,$ see \cite{DR}, proof of Corollary 3.6.

Property $(**)$ for $W$ implies that 
\begin{eqnarray}\label{aass} I^\mu(y)(B_i-1)W^{r}=  I^\mu(y) \left( \begin{array}{c}
                             0 \\
                             \vdots \\
                                                           W \\
                                                          \vdots \\
                              0
                      \end{array}\right) \,\mbox{( $i$-th entry)}.\end{eqnarray}

By assumption on $W,$ one has
\[ I^\mu(y) \left( \begin{array}{c}
                            w_1 \\
                             \vdots \\
                            w_r
                      \end{array}\right)=
      (\sum_{k=1}^r g_k) +g_{\infty}, \]
where $w_1,\ldots,w_r \in W$ and $g_i\in \hat{\K}_i$
(resp. $g_{\infty}  \in \hat{\L}),$
by  Lemma~\ref{Resi}.
Using the monodromy around $t_i$ we get     
\[  I^\mu(y) B_i\left( \begin{array}{c}
                            w_1 \\
                             \vdots \\
                            w_r
                      \end{array}\right)= (\sum_{k\neq i} g_k) +g_{\infty} +
            \lambda g_i. \]
Subtracting theses equalities one obtains 
$ I^\mu(y)(B_i-1)W^{r}\leq \hat{\K}_i$ and (together with 
Equation (\ref{aass}))
$$ I^\mu(y)(B_i-1)W^{r}=I^\mu(y) \left( \begin{array}{c}
                             0 \\
                             \vdots \\
                              W \\
                             \vdots \\
                              0
                      \end{array}\right) \leq \hat{\K}_i.$$

Using the description of $\hat{\K}_i$ in terms of functions,
one sees that for $j=1,\ldots,r,\, j\not= i,$  
$$ I^\mu(y)(B_j-1) \left( \begin{array}{c}
                             0 \\
                             \vdots \\
                              W \\
                             \vdots \\
                              0
                      \end{array}\right)= I^\mu(y) \left( \begin{array}{c}
                             0 \\
                             \vdots \\
                            (A_i-1)  W \\
                             \vdots \\
                              0
                      \end{array}\right) =0,$$
where the expression on the right hand side of the first 
equality is  zero outside the 
$j$-th block entry. Similarly one obtains
$$ I^\mu(y)(B_i-\lambda) \left( \begin{array}{c}
                             0 \\
                             \vdots \\
                              W \\
                             \vdots \\
                              0
                      \end{array}\right)= I^\mu(y) \left( \begin{array}{c}
                             0 \\
                             \vdots \\
                            (A_i-1)  W \\
                             \vdots \\
                              0
                      \end{array}\right) =0.$$

Since $(**)$ holds for $W,$ a block-wise argument 
shows that 
\begin{eqnarray}\label{asss} W^{r} \leq  \ker(I^{\mu}) .\end{eqnarray}

On the other hand,
since $\cap_{i=1}^{r} \ker(A_i|_W-1) = 0$ (Property $(*)$),
we can find an $i\in \{ 1,\ldots,r\}$ and a
solution $f$ in $F\cdot W$ (where 
$F$ is a fundamental system of $D_\a$), such that $f$ has nontrivial 
monodromy at $t_i.$ The Euler transform 
$$g:=I^\mu_{[\alpha_{r+1},\alpha_i]}((\frac{f}{x-t_1},
\ldots,\frac{f}{x-t_r})^{\tr})$$
is a solution of $D_{c_{\mu-1}(\a)}.$ 
 Lemma~\ref{add} implies then that $g$ is not identically zero. 
This gives a contradiction to Equation (\ref{asss}), so $W=0$
and the claim follows.
\Endproof

Finish of the proof of ii):
It follows from the assumptions on the residues and Lemmata 
\ref{Resi} and \ref{hut} that $\dim(\hat{\K})=\dim(\K)$ and 
$\dim(\hat{\L})=\dim(\L).$ It follows then
from Lemma \ref{4.14} that the columns of $I^\mu(y)$ are linearly 
independent.\Endproof

\noindent {\bf Proof} of iii): This follows from dimension 
reasons (using the rank-conditions) 
 and  Lemma \ref{4.14}.\Endproof

\section{Applications of the convolution functors $MC_\lambda$
and $mc_\mu$}\label{apl}

\subsection{Rigid local systems and Fuchsian systems}\label{Rigidesyst}

In this subsection we want to outline 
a construction algorithm for Fuchsian systems
corresponding to irreducible rigid local 
systems under the Riemann-Hilbert correspondence.\\

For 
$\Omega=(\omega_1,\ldots,\omega_r)\in (K^\times)^r,$ the {\em scalar
multiplication with $\Omega$}
$$ \GL_n(K)^r \to  \GL_n(K)^r,\,\, (A_1,\ldots,A_r) \mapsto 
(\omega_1 A_1,\ldots,\omega_rA_r)$$ is denoted by $M_\Omega.$ The 
corresponding 
effect on local systems on the $r$-punctured affine line
is also denoted by $M_\Omega.$
 Similarly,
for $\Delta=(\delta_1,\ldots,\delta_r)\in K^r,$ the {\em 
scalar addition with $\Delta$}
$$ (K^{n\times n})^r \to  (K^{n\times n})^r,\,\, (a_1,
\ldots,a_r) \mapsto (a_1+\delta_1\cdot 1,\ldots,a_r+\delta_r\cdot 1)$$
is denoted by $m_\Delta.$ The corresponding effect on Fuchsian systems is also
denoted by $m_\Delta.$\\

Let $\F$ be a complex irreducible (physically) 
rigid local system. 
By the results of \cite{Katz97}, Chap. 6, and \cite{DR}, Chap. 4,
one can construct $\F$ by applying iteratively 
a suitable sequence of 
scalar multiplications $M_{\Omega^i}$
 (see \cite{DR}, Chap. 4) and middle convolutions
$MC_{\lambda_j}$ 
to a one-dimensional local system $\F_0.$\\

It is easy  to write down a Fuchsian system
$$ D_{\a^0}:Y'= (\frac{a^0_1}{x-t_1}+\cdots + \frac{a^0_r}{x-t_r}) Y, \quad a^0_i \in \CC,$$
 whose solutions 
form the local system $\F_0.$ This system is irreducible and rigid, 
and we (can) assume  that it fulfills the 
assumptions of Theorem \ref{monodromy} iii) 
(i.e., if $a^0_i \in \ZZ$ then $a^0_i=0,$ and 
there exist
at least two elements 
$a^0_{i_1},\, a^0_{i_2}$ such that $a^0_{i_1},\, a^0_{i_2} \notin 
\ZZ$). 
It follows now from Theorem \ref{monodromy} iii)
that there exists
a sequence of scalar additions $$m_{\Delta^i},\,
\Delta^i=(\delta^i_1,\ldots,\delta^i_r),\quad \mbox{\rm such that}\quad
(e^{2\pi i \delta^i_1},
\ldots,e^{2\pi i \delta^i_r}) = \Omega^i,$$
and  middle convolutions
$$mc_{\mu_j},\,\, e^{2\pi i \mu_j}=\lambda_j,$$
such that the iterative application of this sequence to
$D_{\a^0}$ yields an irreducible 
 Fuchsian system $D$ whose monodromy coincides with 
the monodromy of $\F.$ The only thing one has to take care of,
is to choose the scalar additions 
(modulo $\ZZ$)
that the rank 
condition of Theorem \ref{monodromy} iii) is fulfilled.
This is possible in every step by the following argument:

It is shown in \cite{DR}, Chap. 4, resp. Appendix A, that, in the irreducible
case, the 
index of rigidity is preserved by $MC_\lambda,$ resp. 
$mc_\mu$ (which is equal to $2$ in both, the additive and 
the multiplicative, cases). By the compatibility between 
$MC_\lambda$ and $mc_\mu$ (dimension reasons), one can see that if 
two eigenvalues of a matrix which occurs as a component 
in one step of the above ``additive'' construction differ by 
an element of $\ZZ,$ then they correspond to 
a certain Jordan block of length $>1$ in a matrix which
occurs in the  ``multiplicative'' construction, in a way
that the rank condition of Theorem \ref{monodromy} iii) is 
fulfilled.

By the construction 
of $mc_\mu,$ it clear that everything can be 
done in an algorithmic way and is easily 
implemented on the computer. Moreover, one obtains the sections 
of the local 
system $\F$ in a concrete way as iterated integrals,
compare to \cite{HY} and  \cite{Yokoyama02}.\\

\noindent {\bf Remark:} By Crawley-Boevey's solution of the additive 
Deligne-Simpson problem (see \cite{Crawley-Boevey02}) the rigid 
tuples of complex matrices having sum $\equiv 0$ are known and, 
by the additive 
Katz' existence algorithm (see \cite{DR}, Appendix A), these
tuples  can be constructed similar to the above construction.
 
But in general, it is a difficult problem to decide when 
the associated Fuchsian system is irreducible, i.e., the 
associated local system of solutions is an irreducible 
rigid local system. The point is, that in the above construction,
the irreduciblity is ensured by 
Theorem \ref{monodromy} iii), using   the fact that (under the given assumptions)
the functors $MC_{\lambda_j}$ preserve irreduciblity, see 
Theorem \ref{eigen}, iii).

\subsection{Geometric differential
equations}\label{Kurven}

Let 
$X$ be a smooth and geometrically connected algebraic variety
 over an algebraically 
closed field  $K\subseteq \CC,$ 
  $f:Y\to X$ a smooth projective 
morphism and $d$  the universal
differential ${Y}\to \Omega^1_Y.$ 
The 
 Gau\ss-Manin connection 
on  relative 
de Rham cohomology groups $H^i_{\rm DR}(Y/X):=
R^if_*^{\rm DR}({\cal O}_Y,d)$ 
gives rise to a 
 system of differential equations (see \cite{Andre89} for details).
A differential system 
 is said to be {\em arising from 
geometry}
if it is equivalent to an iterated extension
of subfactors of such 
differential systems (see \cite{Andre89}, Chap. II).\\

\begin{thm}\label{coming} Let $K$ be a number field,
 $\a=(a_1,\ldots,a_r),\, a_i\in
 K^{n\times n},\,\mu \in \QQ,$ such that the conditions
of Theorem \ref{monodromy} hold for $D_\a.$ 
If $D_\a$ is arising from geometry, then $D_{c_\mu(\a)}$
(resp. $D_{mc_\mu(\a)}$) is arising from geometry. 
\end{thm}

\proof
The claim  follows from the construction 
of the period matrix and the result of Andr\'e \cite{Andre89}, 
saying that the category of differential modules which arise from 
geometry is closed under taking higher direct images.\Endproof \\

Let us consider an example:

\begin{lemma}
Let $$p(x)=4 (x-t_1)(x-t_2)(x-t_3)=4x^3-g_2x-g_3,$$ $B\in \CC$ and 
 $$ L_n:=L_n(p,B):= p(x)y{''} + \frac{1}{2} {p'(x)}y'
    -(n (n+1)x+ B)y$$
 the Lam\'e differential equation of index $n \in \QQ.$ 
 Then, $L_n$ can be transformed into the Fuchsian system

\begin{eqnarray}Y'& =&\sum_{i=1}^3 {a_i \over x-t_i} Y\nonumber \\
&:=&
     \left({1 \over{x-t_1}} \left( \begin{array}{cc}
                            0 & 1 \\
                             0 & \frac{1}{2} 
                      \end{array} \right)
     +
      {1 \over{x-t_2}} \left( \begin{array}{cc}
                               0 & 0 \\
                               l_1 & -\frac{1}{2} \\
                              
                      \end{array} \right)+ 
         {1 \over{x-t_3}} \left( \begin{array}{cc}
                               0 & 0 \\
                            l_2 & -\frac{1}{2} \\
                      \end{array} \right)\right)Y, \nonumber\end{eqnarray}
where
 \[l_1={ t_2 n (n+1)+ B \over 4 (t_2-t_3)}
\quad \mbox{and} \quad l_1+l_2= \frac{n (n+1)}{4}. \]
\end{lemma}

\proof
 The differential system which corresponds to 
 $L_n$ is
 \[ Y'=     \left( \begin{array}{cc}
                            0 & 1 \\
              { n (n+1)x+B \over p(x)}         & -\frac{1}{2}{p'(x) \over{p(x)}} 
                      \end{array} \right)Y .\]
 Using the gauge transformation $Y \mapsto 
   \left( \begin{array}{cc}
                            1 & 0 \\
                             0 & x-t_1 
                      \end{array} \right)Y $
we get the equivalent system
 \[ Y'=   \left( {1 \over{x-t_1}} \left( \begin{array}{cc}
                            0 & 1 \\
                             0 & \frac{1}{2} 
                      \end{array} \right)
             + \left( \begin{array}{cc}
                            0 & 0 \\
              { n (n+1)x+ B \over 4 (x-t_2)(x-t_3)}  & -\frac{1}{2}\sum_{i=2}^3
                {1\over{x-t_i}} 
                      \end{array} \right)\right)Y. \]
Since
  \begin{eqnarray}  { n (n+1)x+ B \over 4 (x-t_2)(x-t_3)} &=&
   { t_2 n (n+1)+ B   \over 4
         (x-t_2)(t_2-t_3)} - { t_3 n (n+1)+ B  \over 4
         (x-t_3)(t_2-t_3)}, \nonumber \end{eqnarray}
the claim follows.
\Endproof

\begin{defn}{\rm  We say that a  system  of differential equations $D$
is 
in {\em Okubo normal form}, if
$$ D:Y'=(x-T)^{-1}bY \quad \mbox{\rm 
or, equivalently,}\quad  D:(x-T)Y'=bY,$$ 
where $b\in  \CC^{n\times n}$ and $T$ is a 
diagonal matrix 
 $T={\rm diag}(t_1,\ldots,t_n),\, t_i \in \CC$
(here possibly $t_i=t_j$ for $i\not= j$).
}\end{defn}

\begin{lemma}\label{lame} Let $r\ge 3,$ $\mu\in \CC \setminus \ZZ$ and 
$\a:=(a_1,\ldots,a_r),\, a_i \in \CC^{2\times 2},$ 
where $a_1,a_2,a_3$ are as in the previous lemma. If   $rk(a_i)=2$ for 
$i>3,$ and if 
 $-\mu $ is no eigenvalue of 
$ a_1+\cdots +a_r,$ then $D_{mc_\mu(\a)}$
is equivalent to  the following differential system in Okubo form:
 \[ D(L_n,\a,\mu):(x-T)Y'=(\tilde{c}+\mu)Y\]
where $T={\rm diag}(t_1,t_2,t_3,t_4,t_4,\ldots,t_r,t_r)$ and 
\[\tilde{c}= \left( \begin{array}{ccc|c|c|c}
                            {1\over 2} & -{1\over 2} & -{1\over 2} & 
                        \begin{array}{cc} (0,1)a_4 \end{array}  &   & 
                         \begin{array}{cc}   (0,1)a_r \end{array}\\
                             -2l_1 +{1\over 2} & -{1\over 2} & -{1\over 2} & 
                            \begin{array}{cc} (-2l_1,1)a_4 \end{array}  & \ldots &
                            \begin{array}{cc} (-2l_1,1)a_r\end{array}  \\
                              -2l_2 +{1\over 2} & -{1\over 2} & -{1\over 2} & \begin{array}{cc} (-2l_2,1)a_4
 \end{array}& &\begin{array}{cc}(-2l_2,1)a_r \end{array}\\
                             &&&&&\\  \hline &&&&&\\
                       \begin{array}{c}
                                 1 \\
                              {1\over 2} 
                            \end{array} &
                       \begin{array}{c}
                                 0 \\
                             - {1\over 2} 
                            \end{array} &
                      \begin{array}{c}
                                 0 \\
                             - {1\over 2} 
                            \end{array} &
                         a_4  & \ldots & a_r \\
                             &&&&&\\  \hline &&&&&\\ 
                                  \vdots & \vdots &\vdots & \vdots  & &\vdots\\
                       &&&&&\\  \hline &&&&&\\ 
                        \begin{array}{c}
                                 1 \\
                              {1\over 2} 
                            \end{array} &
                       \begin{array}{c}
                                 0 \\
                             - {1\over 2} 
                            \end{array} &
                      \begin{array}{c}
                                 0 \\
                             - {1\over 2} 
                            \end{array} &
                         a_4  & \ldots & a_r \\
                      \end{array} \right) \]
\end{lemma}

\proof In the notation of Subsection \ref{addi}. 
Let $c_0(\a)=(b_1,\ldots,b_r),$ where $b_i\in \CC^{2r\times 2r}.$ 
One has $D_{c_\mu(\a)}:Y'=(x-T)^{-1}(b+\mu)Y,$ where 
$b:=b_1+\ldots+ b_r.$  

 Let
\[   B_1:=  \left( \begin{array}{cc}
                            1 & -2 \\
                             0 & 1 
                      \end{array} \right),
    B_2:=  \left( \begin{array}{cc}
                            1 & 0 \\
                            -2l_1 & 1 
                      \end{array} \right),
   B_3:=  \left( \begin{array}{cc}
                            1 & 0 \\
                            -2l_2 & 1 
                      \end{array} \right). \]

Then
\[ B_1 a_1 B_1^{-1}=  \left( \begin{array}{cc}
                            0 & 0 \\
                            0 & \frac{1}{2} 
                      \end{array} \right)\quad {\rm and } \quad B_2 a_2 B_2^{-1}=  B_3 a_3 B_3^{-1}= \left( \begin{array}{cc}
                            0 & 0 \\
                            0 & -\frac{1}{2} 
                      \end{array} \right). \]
Consider the block-diagonalmatrix $d:={\rm diag}(B_1,B_2,B_3,E_2,\ldots,E_2) .$
One computes that 
\[dbd^{-1}= \left( \begin{array}{ccc|ccccccc}
      B_1 a_1B_1^{-1} & B_1 a_2 B_2^{-1} & B_1 a_3B_3^{-1}& B_1 a_4 & \ldots&
      B_1 a_r 
      \\
   B_2 a_1B_1^{-1} & B_2 a_2 B_2^{-1} & B_2 a_3B_3^{-1}& B_2 a_4 &
     \ldots& B_2 a_r   \\
     B_3 a_1B_1^{-1} & B_3 a_2 B_2^{-1} & B_3 a_3B_3^{-1}& B_3 a_4 & 
        \ldots& B_3 a_r  \\
            \hline
     a_1  B_1^{-1} &  a_2 B_2^{-1} &  a_3 B_3^{-1}&  a_4 &  \dots & a_r   \\
          \vdots         &   \vdots &  \vdots & &      \vdots & & \\
       a_1  B_1^{-1} &  a_2 B_2^{-1} &  a_3 B_3^{-1}&  a_4 &   \dots & a_r \\
     
                      \end{array} \right) \]
It is easily checked that conjugating 
$b$ with $d$ gives an equivalence between 
$D_{c_\mu(\a)}$ and $Y'=(x-T)^{-1}(c+\mu)Y$
with $c=dbd^{-1}.$

Under 
the assumptions, the space ${\frak l}\leq \CC^{2r}$ is zero. 
Since $B_ia_iB_i^{-1}$
is diagonal, 
factoring out the space ${\frak k}$ corresponds to canceling the
first, third and fifth row and column of $c.$ This yields  $\tilde{c}.$
\Endproof

Baldassarri \cite{Baldassari80} gives examples of  
Lam\'e equations with finite monodromy. 
E.g., it is shown that  the monodromy group
of  $L_n(p(x),B),$ ($n=1,$ $p(x)=4x^3+g_3,$ $B=0$), is the 
symmetric group on $3$ letters. 
Beukers and  van der Waall \cite{BW} 
list all finite groups which can occur as the 
monodromy groups of  Lam\'e equations and 
 they give many examples of such equations. 
Also,  van der Waall \cite{vdW} gives an 
algorithm  to detect all Lam\'e equations with finite monodromy group.
One equation
 which can be found 
in \cite{BW}, Table 4, is given 
by $L_n(p(x),B)$ ($n=1/6,$ $p(x)=4x^3-x,$ $B=0$), whose monodromy 
group is isomorphic to the complex reflection group 
$G_{13}.$

 This yields a new family of differential systems 
arising from geometry:

\begin{cor}\label{aka} Under the assumptions of Lemma \ref{lame}.
Let  
  $L_n(p,B)$ a Lam\'e equation which 
has finite monodromy (e.g., if $n=1,$ $p(x)=4x^3+g_3,$ $B=0,$ or
$n=1/6,$ $p(x)=4x^3-x,$ $B=0$), $a_4,\ldots,a_r$ scalar matrices  
contained in $\QQ^{2\times 2}$ 
 and $\mu \in \QQ.$ Then the following holds:\\

i) The differential 
system $D(L_n(p,B),\a,\mu)$ is arising from geometry.\\

ii) The solutions of $D(L_n(p,B),\a,\mu)$ are $G$-functions.\\
\end{cor}

\proof This follows from Theorem \ref{coming} and 
\cite{Andre89}, Chap. V. \Endproof

\subsection{Transformation of the $p$-curvature under
$mc_\mu$}\label{pcurv}

In this section, 
we want to study how the $p$-curvature changes under the convolution:\\

Let $K$ be a number field and 
$D:Y'=aY,$ where $a=(a_{i,j})\in K(x)^{n\times n}.$
Successive application of
differentiation yields differential systems 
$$D^{(n)}: Y^{(n)}=\widehat{a(n)}Y.$$
In the following, $\p$ always denotes a prime 
of $K$ which lies over $p.$
For almost all primes $p,$   one can 
reduce $\widehat{a(p)}$ modulo 
$\p,$ 
 in order 
 to obtain the {\em $p$-curvature matrices}
$$a(\p):=\widehat{a(p)} \mod \p.$$
The $p$-curvature matrices
encode many  arithmetic properties of the differential system 
 $D$ and are 
conjecturally related to questions about the  
geometric nature of $D:$

\begin{con}\label{Gro} i) {\rm (Grothendieck-Katz, see \cite{Katz82}, \cite{Andre96})}
The Lie algebra of the differential Galois group of $D$ is minimal
to the property that, for almost all primes $\p$
of $K,$ its reduction modulo $\p$ contains 
the $p$-curvature matrix $a(\p).$

ii) {\rm (Bombieri-Dwork, see \cite{Andre89})} If $D$ 
is globally nilpotent, i.e., ${a}{(\p)}$ is nilpotent
 for almost all primes $\p,$ then $D$ is arising from geometry.
\end{con}

\noindent {\bf Remark:} (i) 
The Grothendieck-Katz conjecture
implies the $p$-curvature conjecture of Grothendieck: If 
$a(\p)=0$ for almost all $\p,$ then $D$ has a fundamental set of solutions 
consisting of algebraic functions.

(ii) It is well known that if 
$D$ has a fundamental set of solutions consisting of 
algebraic functions,
then ${a}{(\p)}= 0$ for almost all $\p.$ Also, if $D$ is
arising from geometry, then ${a}{(\p)}$ is nilpotent 
for almost all $\p.$\\

\begin{rem} {\rm (Okubo)}
 Let $ D:(x-T)Y'=bY$
be   a system of differential equations in Okubo normal form. 
 Then, $(x-T)Y^{(2)}=(b-1)Y'.$
\end{rem}

An induction  yields the following recursion formula for 
the $p$-curvature matrix of a system of differential equations 
in Okubo normal:

\begin{lemma}\label{recu} Let $ D:(x-T)Y'=bY$
be   a system of differential equations in Okubo normal form. 
 Then
 $$\widehat{a(n)}=(x-T)^{-1}(b-n+1)\cdot (x-T)^{-1}(b-n+2) \cdots (x-T)^{-1}
(b-1)\cdot (x-T)^{-1}b.$$
\end{lemma}

\begin{thm}\label{nilpo} Let  
$\a=(a_1,\ldots,a_r), \, a_i \in K^{n\times n},$ such that 
 $$a(\p)^k= 0.$$ Let $\mu \in \QQ$ 
and denote by $c_\mu(a(\p))$
(resp. $mc_\mu(a(\p))$) the $p$-curvature matrix
of $D_{c_\mu(\a)}$ (resp. $D_{mc_\mu(\a)}$). \\

i) If $\mu=-1,$ then $c_\mu(a(\p))^{k+1}=0$ and 
 $mc_\mu(a(\p))^{k+1}= 0.$\\

ii) If $\mu =\frac{n_1}{n_2}$ and $p$ does not divide $n_1n_2,$
 then   $$c_{\mu-1}(a(\p))^{k+2}= 0\,\,\,\, and \,\,\,\,
 mc_{\mu-1}(a(\p))^{k+2}=0.$$

\end{thm}

\proof The convolution of $D_\a$
is a differential system in Okubo normal form:
$$D_{c_\mu(\a)}:Y'=\sum_{k=1}^r \frac{b_k}{x-t_k}Y
=(x-T)^{-1}(\sum_{k=1}^r b_k)Y,$$
where $T$ is the diagonal matrix 
$T={\rm diag}(t_1,\ldots,t_1,\ldots,t_r,\ldots,t_r)$
(every $t_k$ occurs $n$ times) and $b_k$ is as in 
Section \ref{addi}. 
If $\mu=-1$ then $(\sum_{k=1}^r b_k)$ is a blockmatrix $b=(b_{i,j})$
with $b_{i,j}=a_j -\delta_{i,j}E_n.$

 Using the gauge transformations with $(x-T)$ and
 \[H:=\left( \begin{array}{cccc}
       E_n & -E_n & 0 & \ldots \\
       0  &      \ddots&\ddots    \\
       \vdots&  & &- E_n \\
       0  & \ldots&& E_n
    \end{array} \right) \]
one sees that $D_{c_{-1}(\a)}$ is equivalent to the following 
system:
 \[   
  Y'= \left( \begin{array}{cccc}
       0 &  \ldots & \ldots & 0 \\
            \vdots &&&\vdots\\
                         0 & \ldots & \ldots &0 \\
      {{a_1}\over{x-t_1}} &({{a_1}\over{x-t_1}}+{{a_2}\over{x-t_2}}) 
&\ldots & ({a_{1}\over{x-t_1}}+\cdots + 
{a_{r}\over{x-t_r}})
    \end{array} \right)Y. \]
Thus $c_{-1}(a(\p)) $ is equivalent to
 \[  \left( \begin{array}{cccc}
       0 &  \ldots & \ldots &0 \\
            \vdots &&&\vdots\\
                         0 & \ldots&\ldots  & 0 \\
      \ast & \ast &\ast  & a(\p)
    \end{array} \right) .  \]
It follows from $a(\p)^k= 0$ that 
 $c_{-1}(a(\p))^{k+1} = 0 $ and i) follows.

Let $\hat{\mu}\in {\Bbb N}^+$ be the smallest natural number 
such that $\hat{\mu}\equiv \mu \mod \p$ and 
$b_\infty:=b_1+\cdots+b_r-\mu,$ where $c_\mu(\a)=(b_1,\ldots,b_r).$ 
 Let further 
 $$h_1:=(x-T)^{-1} b_\infty\cdot (x-T)^{-1} (b_\infty+1) \cdots  (x-T)^{-1} (b_\infty-1+\hat{\mu})$$ and
$$h_2:=(x-T)^{-1} (b_\infty+\hat{\mu})\cdot (x-T)^{-1} (b_\infty+\hat{\mu}+1) \cdots  (x-T)^{-1} (b_\infty-1+p).$$
Then $c_{\mu-1}(a(\p))=h_2h_1\mod \p$ and $c_{-1}(a(\p))=(h_1h_2\,\mod \p)$
by the above recursion formula (Lemma \ref{recu}).
It follows from $$c_{-1}(a(\p))^{k+1}=((h_1h_2)^{k+1}\,\mod \p)= 0,$$
that 
$$c_{\mu -1}(a(\p))^{k+2}=((h_2h_1)^{k+2}\,\mod \p)= 0,$$
giving ii). 
\Endproof

\begin{cor} Let  $\a=(a_1,\ldots,a_r),\, a_i\in K^{n\times n},$ such that 
$D_\a$ is a globally nilpotent  and $\mu \in \QQ.$
Then  $D_{c_\mu(\a)}$ (resp. $D_{mc_\mu(\a)}$)
is globally nilpotent. 
\end{cor}

\begin{lemma}\label{ak} Let $D(L_n(p(x),B),\a,\mu)$
be as in Corollary \ref{aka}. \\

i) The Grothendieck-Katz conjecture 
is true for $D(L_n,\a,\mu).$\\

ii) The system $D(L_n,\a,\mu)$ is globally nilpotent of rank 3.\\
\end{lemma}

\proof i) In the notation of Section \ref{s2}.
Let 
$G$ be the monodromy group of $D_\a.$ Let 
$X=\CC\setminus T$ (remember that
$t_1,t_2,t_3$ are determined by $L_n$),
$X_1 \to X$ be the unramified cover of $X$ which is associated 
to the homomorphism $\pi_1(X)\to G \leq \GL_2(\CC)$ and 
$X_2$ the cyclic cover of $\CC^\times$ which is associated to 
$\pi_1(\CC^\times)\to \CC^\times,\, \gamma \mapsto e^{2\pi i \mu}.$
 Let $Y_1:=X_1\times_X E ,$
 $Y_2:=X_2\times_{\CC^\times} E$ and 
$\tilde{Y}:=Y_1\times_E Y_2.$  By construction,
 $\tilde{Y}$ is an unramified 
cover of $E$ and admits, via ${\rm p}_2,$ a map $\tilde{f}:\tilde{Y}\to X.$
Let $Y$ denote the compactification of $\tilde{Y}$
with respect to the first coordinate and $f:Y\to X$ the morphism induced
by $\tilde{f}.$ It follows from Riemann's existence theorem that
$f$ arises
 from  an underlying  smooth  map of 
varieties (i.e., $f$ is the effect on the complex points), 
denoted by $\hat{f}:\hat{Y}\to \hat{X},$ where $\hat{X}$ and $\hat{Y}$ are
smooth connected varieties  over some  number field $K.$
 It follows from the 
Leray spectral sequence that
$D_{mc_{\mu}(\a)}$ is a differential system 
which is equivalent to a subfactor of the Gau\ss-Manin 
connection $\nabla: H^1_{\rm DR}(\hat{Y}/\hat{X})
\to \Omega^1_{\hat{X}}\otimes H^1_{\rm DR}(\hat{Y}/\hat{X}).$

By \cite{Andre02}, Theorem 0.7.1, the claim follows from the 
connectivity of the motivic Galois group of at least one 
fibre $\hat{Y}_s$
(which is a nonsingular curve in our case),
 where $s$ is a geometric point  of $\hat{X}.$
 But this follows analogously to \cite{Andre02}, Ex. 16.3,
from the results of $\cite{Andre96},$ 
relating the motivic Galois group of $\hat{Y}_s$ to the Mumford-Tate
group of the Jacobian of $\hat{Y}_s.$

The claim ii) follows from Theorem \ref{nilpo}, ii).
\Endproof

\bibliographystyle{plain} \bibliography{p}

Michael Dettweiler 

IWR, Universit\"at Heidelberg, 

INF 368

69121 Heidelberg, Deutschland

e-mail: michael.dettweiler@iwr.uni-heidelberg.de\\

Stefan Reiter

G.R.I.M.M., Universit\'e de Toulouse II 

5. All\'ees de A. Machado

31058 Toulouse, France

email: reiter@math.jussieu.fr

\end{document}